\documentclass[11pt,a4paper]{amsart}
\usepackage{amsmath,amsfonts,amsthm,amsopn,color,amssymb,enumitem}
\usepackage[a4paper]{geometry}
\usepackage{palatino}
\usepackage{graphicx}
\usepackage[colorlinks=true]{hyperref}
\hypersetup{urlcolor=blue, citecolor=red, linkcolor=blue}

\usepackage{cite}
\usepackage{relsize}
\usepackage{esint}
\usepackage{verbatim}
\usepackage{mathrsfs}
\usepackage{xcolor}
\usepackage{tikz}

\usepackage[leqno]{amsmath}
\usepackage{times}
\usepackage{cancel}
\usepackage[colorinlistoftodos]{todonotes}

\newtheorem{Thm}{Theorem}[section]

\newtheorem{Lem}[Thm]{Lemma}

\newtheorem{Prop}[Thm]{Proposition}

\newtheorem{Rem}[Thm]{Remark}

\numberwithin{equation}{section}

\newcommand{{\X}}{{\widetilde X}}
\newcommand{{\Y}}{{\widetilde Y}}

\newcommand{\R}{\mathbb{R}}

\newcommand{\cA}{{\mathcal A}}

\newcommand{\cD}{{\mathcal D}}

\newcommand{\cF}{{\mathcal F}}

\newcommand{\cH}{{\mathcal H}}

\newcommand{\cJ}{{\mathcal J}}
\newcommand{\cK}{{\mathcal K}}
\newcommand{\cL}{{\mathcal L}}

\newcommand{\cN}{{\mathcal N}}

\newcommand{\cR}{{\mathcal R}}
\newcommand{\cS}{{\mathcal S}}

\newcommand{\al}{\alpha}
\newcommand{\be}{\beta}
\newcommand{\ga}{\gamma}
\newcommand{\de}{\delta}
\newcommand{\ka}{\kappa}
\newcommand{\la}{\lambda}

\newcommand{\De}{\Delta}
\newcommand{\Ga}{\Gamma}
\newcommand{\Om}{\Omega}
\newcommand{\om}{\omega}
\newcommand{\Si}{\Sigma}

\newcommand{\eps}{\varepsilon}
\newcommand{\pa}{\partial}

\newcommand{\ubar}{\bar u}
\newcommand{\dd}{\,\mathrm d}

\newcommand{\wt}{\widetilde}

\newcommand{\beq[1]}{\begin{equation}\label{eq:#1}}
\newcommand{\eeq}{\end{equation}}

\newenvironment{altproof}[1]
{\noindent
{\em Proof of {#1}}.}
{\nopagebreak\mbox{}\hfill $\Box$\par\addvspace{0.5cm}}

\begin{document}

\title{Periodic vortex dipoles near closed geodesics on Riemannian surfaces}
\author[T. Bartsch]{Thomas Bartsch$^\dagger$}\thanks{$\dagger$Supported by DFG grant BA 1009/19-1. T.B. thanks the university Roma La Sapienza for support and hospitality.}
\address{\noindent Mathematisches Institut, Universit\"at Giessen}
\email{Thomas.Bartsch@math.uni-giessen.de}
\author[A. Pistoia]{Angela Pistoia$^\ddagger$}\thanks{$\ddagger$Supported by INDAM-GNAMPA project 2025}
\address{\noindent  Dipartimento di Scienze di Base e Applicate per l'Ingegneria, Università degli Studi di Roma Sapienza}
\email{angela.pistoia@uniroma1.it}

\date{}
\begin{abstract}
We study the motion of a pair $(p_1,p_2)\in\Si\times\Si$ of point vortices with equal magnitude and opposite signs on a closed oriented Riemannian surface $(\Si,g)$. Given a closed geodesic that satisfies a non-degeneracy condition we construct a one-parameter family of periodic solutions $(p_1^\eps,p_2^\eps)$, $0<\eps<\eps_0$, of the Hamiltonian  dynamical system for point vortex pairs whose trajectories converge towards the geodesic curve. The two vortices lie on opposite sides of the geodesic, are separated of order $\eps$ and have period $2\pi\eps$. The proof is based on a differentiated expansion of the Riemannian distance in Fermi coordinates, an augmented Lyapunov--Schmidt reduction, and an intricate contraction argument in scaled variables. The reduced normal equation is governed by the Jacobi operator of the closed geodesic, and the first order correction to the normal center of the dipole is characterized by an inhomogeneous Jacobi equation involving the curvature and the regular part of the Green function of the Laplace-Beltrami operator.
\end{abstract}

\keywords{vortex dynamics, point vortices, vortex dipoles, periodic solutions, Riemannian surfaces, Hamiltonian systems, Jacobi operator, Synge world function}
\subjclass[2020]{Primary 76B47; Secondary 37J46, 53C22}


\maketitle

%
\section{Introduction}\label{sec:intro}
Point-vortex dynamics provides a finite-dimensional Hamiltonian description of
the interaction of concentrated vortices in an ideal incompressible fluid. It has its roots in the work of Helmholtz~\cite{Helmholtz:1858} and Kirchhoff~\cite{Kirchhoff:1876}, and has been studied by many mathematicians and physicists. Most papers deal with point vortex dynamics 
in the plane or on the round sphere, motivated by various applications, from thin films to geophysics and astrophysics. General accounts of 
point-vortex dynamics can be found in
\cite{newton:2001,saffman:1992}. Geometric formulations on curved surfaces were developed, in particular, 
by Boatto and Koiller \cite{Boatto-Koiller:2015}, Dritschel and Boatto
\cite{Dritschel-Boatto:2015}, and Gustafsson \cite{Gustafsson:2022}.
On a Riemannian surface point vortex dynamics reflects both the local geometry of the
metric and the global structure of the Green function of the Laplace--Beltrami operator.
A particularly interesting regime arises when two vortices of equal magnitude and
opposite signs approach one another. Then the singular logarithmic interaction 
produces a rapidly travelling dipole, and its leading-order motion is governed
by the geodesic flow. The purpose of this paper is to construct 
periodic dipoles near prescribed nondegenerate closed geodesics.

Let $(\Si,g)$ be a closed oriented Riemannian surface and let $\cF_2\Si:=\{(p_1,p_2)\in\Si\times\Si:p_1\ne p_2\}$ be the configuration space of two different points in $\Si$. We consider the Green--Robin point-vortex model, in which the dynamics of two point vortices $p=(p_1,p_2)\in\cF_2\Si$ with vortex strengths $\ka_1,\ka_2\in\R$ is governed by a Hamiltonian system
\begin{equation}\label{HS}
 \ka\dot p = X_\cH(p)
\end{equation}
with $\ka=\begin{pmatrix}\ka_1 &0\\ 0 &\ka_2 \end{pmatrix}$ and Hamiltonian $\cH:\cF_2\Si\to\R$ given by
\begin{equation}\label{def:H}
  \cH(p_1,p_2) = -2\ka_1\ka_2 G_\Si(p_1,p_2)  + \ka_1^2R_\Si(p_1) + \ka_2^2R_\Si(p_2).
\end{equation}
Here $G_\Si$ denotes the Green function for the Laplace--Beltrami operator on $\Si$, and $R_\Si$ the Robin function, i.e.\ with $d_g:\Si\times\Si\to\R$ denoting the distance function: 
\[
  R_\Si(x) = \lim_{y\to x}\left(-\frac{1}{4\pi}\ln d_g^2(x,y) - G_\Si(x,y)\right).
\]
In order to define the vector field $X_\cH$ let $\Om$ be the area form on $(\Si,g)$, and let $\om := \pi_1^*(\Om)+\pi_2^*(\Om)$ be the induced symplectic form on $\Si\times\Si$; here $\pi_i:\Si\times\Si\to\Si$ is the projection onto the $i$-th factor. Then $X_\cH$ is defined by the equation
\[
  \om\big(X_\cH(p),\cdot\,\big) = -d\cH(p).
\]

The relation between geodesics and dipoles, where $\ka_1=-\ka_2$ and $p_1$, $p_2$ are close, goes back to Kimura
\cite{Kimura:1999}, who studied vortex motion on surfaces of constant
curvature and conjectured that an infinitesimally close pair of opposite
vortices follows a geodesic also on a surface of variable curvature; a vortex 
dipole may be considered a "geometry checker". Since the speed diverges as the distance of the vortices tends to zero,
the limiting dynamics must be described after a suitable rescaling of time.
Koiller and Boatto \cite{Koiller-Boatto:2009} developed a geometric approach to this limit.
Gustafsson \cite{Gustafsson:2022}, in particular Theorem~8.1 therein,
derived the geodesic equation directly from the dipole-limit equations in
local complex coordinates. More recently, Drivas, Glukhovskiy and Khesin
\cite{Drivas-Glukhovskiy-Khesin:2024} studied vortex pairs for general vorticities
$\ka_1,\ka_2$ and obtained magnetic geodesics as the limiting
motion when the distance of the vortices tends to $0$. In the dipole case $\ka_1=-\ka_2$ the magnetic term disappears,
yielding another proof of the geodesic limit predicted by Kimura.

The dynamics for small but nonzero separation has also been investigated
in specific geometries. Regis Rodrigues, Castilho and Koiller
\cite{Rodrigues-Castilho-Koiller:2018} studied opposite vortex pairs on
a triaxial ellipsoid and numerically tested their relation to the geodesic
flow. Sakajo and Shimizu \cite{Sakajo-Shimizu:2016} analyzed the two-vortex
problem on a torus, including explicit symmetric dipole configurations
with periodic motion and their geodesic limits. Thus periodic dipoles
are explicitly known for particular symmetric surfaces. Periodic solutions of the Hamiltonian system for 
$N$ point vortices on a Riemannian sphere equipped with a generic metric have been conctructed by 
Wang \cite{wang:2021} provided all vorticities are positive and satisfy a certain commensurability condition. 
The existence of periodic solutions for $N\ge2$ 
point vortices has also been studied in bounded planar domains by Bartsch and Dai \cite{bartsch-dai:2016}, 
Bartsch and Gebhard \cite{bartsch-gebhard:2017}, Bartsch and Sacchet \cite{bartsch-sacchet:2018}. 
In these papers the authors used singular perturbation methods to construct families of periodic solutions where the distance of the vortices tends to zero. However, a basic assumption was that the total vorticity is not $0$, hence the case of two vortices with total vorticity $\ka_1+\ka_2=0$ was explicitly excluded.

The geodesic limit in the case $\ka_1=-\ka_2$ for two point vortices suggests a natural persistence problem: does a given
closed geodesic give rise to periodic orbits of the vortex-pair system for initial values close to the geodesic? Convergence to the geodesic flow 
does not settle this question. Even if the limiting trajectory is closed, the perturbed trajectory need not be closed, and may even leave a 
neighborhood of the geodesic. Exact periodicity requires solving \eqref{HS} with periodic boundary conditions, for initial values approaching a point on the closed geodesic. 

Our main result, Theorem~\ref{thm:main}, answers this question affirmatively for every closed geodesic $\Ga$ whose normal Jacobi
operator is invertible. This result thus relates periodic solutions of the second order differential equation for the geodesic to those of the first order Hamiltonian system \eqref{HS}. More precisely, if $\ga$ is an arc-length parametrization of $\Ga$, the nondegeneracy assumption is that
\[
 h''+K(\ga(x))h=0
\]
has no nontrivial periodic solution having the same period as the arc-length parametrization of $\Ga$; here $K$ denotes the Gauss curvature.
After normalizing the length of $\Ga$ to one and the vortex strengths to
$\ka_1=1=-\ka_2$, we construct a smooth family of periodic solutions
with period $2\pi\eps$, for every sufficiently small $\eps>0$.
In Fermi coordinates $(x,y)$ around $\Ga$, the vortices travel on opposite
sides of the geodesic and satisfy
\[
 x_i^\eps(t)=\frac{t}{2\pi\eps}+O(\eps^2),\qquad
 y_1^\eps(t)=\eps+O(\eps^2),\qquad
 y_2^\eps(t)=-\eps+O(\eps^2).
\]
The estimates are uniform in time, and the theorem gives the corresponding
Sobolev estimates after rescaling to a fixed period. In particular, both
vortices converge to the same closed geodesic under this rescaling.
The construction requires neither a symmetry of the surface nor a sign
condition on its curvature. Its geometric hypothesis is precisely the
absence of periodic normal Jacobi fields along the chosen geodesic.

Beyond existence, the construction describes the first corrections to
the limiting motion. The tangential separation $|x_1^\eps-x_2^\eps|$ is of order $\eps^3$,
whereas the normal separation $|y_1^\eps-y_2^\eps|$ is of order $2\eps+O(\eps^3)$. The leading term of the normal-center corrections
$|y_1^\eps-\eps|$, $|y_2^\eps+\eps|$ is of order $\eps^2$ and is determined by an inhomogeneous Jacobi equation; see
\eqref{eq:Zsolution}--\eqref{eq:beta-limit-final}. Its forcing involves
the normal derivative of the curvature and derivatives of the regular
part of the Hamiltonian. This identifies where the global Green function
enters the correction to the locally determined geodesic motion.

The proof combines a differentiated expansion up to order 4 of the squared Riemannian
distance $d_g^2(p_1,p_2)$ in Fermi coordinates with an augmented Lyapunov-Schmidt type reduction of the periodic
system, and an intricate contraction mapping argument in sutiably scaled variables. The function $\frac12d_g^2$ is Synge's world function in general relativity. Although it has been investigated in the literature, 
see \cite{Avramidi:2015} for instance, we did not find our expansion including bounds for the remainder terms appearing in the expansion of the derivatives of $d_g^2$.

The paper is organized as follows. In the next section we state and discuss our main result. There we also give a sketch of the proof. In Section~\ref{sec:prelim} we recall some facts about the Green function of the Laplace-Beltrami operator on $\Si$ and formulate system \eqref{HS} in terms of Fermi coordinates. An important ingredient in the proof  of the main result is the expansion of $d_g^2$ and its derivatives up to order $4$ in Fermi coordinates near the geodesic. This will be presented in Section~\ref{sec:Taylor}. Finally Section~\ref{sec:proof} contains the proof of our main result.

Throughout the paper we assume without loss of generality $\ka_1=-\ka_2=1$ in \eqref{HS}.

\section{Statement of results}\label{sec:results}
Let $\Ga$ be a closed geodesic of length one on an oriented Riemannian surface $(\Si,g)$, and let $\ga:\R\to\Ga\subset\Si$ be its unit-speed parametrization. Choose the unit normal field $\nu$ along $\Ga$ so that $\Om(\dot\ga,\nu)=1$, where $\Om$ is the area form. Then $\ga$ and $\nu$ are $1$-periodic. We do not assume that $\Ga$ is embedded, in particular $\ga|_{[0,1)}$ need not be injective. 

Let us first recall the Fermi coordinates near $\Ga$. These are defined by
\[
  \Phi:\cS := \R\times(-\de,\de)\to\Si,\quad \Phi(x,y) := \exp_{\ga(x)}(y\nu(x)),
\]
for some $\de>0$. Here ${\exp}_{\ga(x)}:T_{\ga(x)}\Si\to \Si$ is the exponential map in $\Si$ at the point $\ga(x)\in\Si$. We remark that
$\Phi(x,0)=\ga(x)$, hence $\Phi(x+1,0)=\Phi(x,0)$. 

Let $K:\Si\to\R$ denote the Gauss curvature of $\Si$ and set $\ka(x):=K\big(\ga(x)\big)$ for $x\in\R$. Let $H^1_{per}$ denote the Hilbert space of $1$-periodic functions in $H^1_{loc}(\R,\R)$, and define the operator
\begin{equation}\label{def:J_Ga}
  \cJ_\Ga: H^1_{per}\to H^{-1}_{per}, \quad\cJ_\Ga(h) := -h'' - \ka h.
\end{equation}
We will also write $H^1_{2\pi}$ for the Hilbert space of $2\pi$-periodic functions in $H^1_{loc}(\R,\R)${\color{red}.} Now we can state our main result.

\begin{Thm}\label{thm:main}
Let $(\Sigma,g)$ be a closed oriented Riemannian surface and let $\Ga$ be a closed
geodesic of length one. Assume that $\Ga$ is nondegenerate in the sense that $\ker \cJ_\Ga=\{0\}$. Then there exists $\eps_0>0$ and a smooth map 
\[
  (0,\eps_0) \to \big(H^1_{2\pi}\big)^4\times\R,\quad \eps\mapsto \big(\phi_1^\eps,\phi_2^\eps,\psi_1^\eps,\psi_2^\eps,\be_\eps\big)
\] 
such that, for every $0<\eps<\eps_0$, the vortex-pair Hamiltonian system \eqref{HS} with $\ka_1=1=-\ka_2$ has a periodic solution $p^\eps(t)=\big(\Phi\big(z_1^\eps(t)\big),\Phi\big(z_2^\eps(t)\big)\big)$
of the form $z_i^\eps=(x_i^\eps,y_i^\eps):\R\to\cS$ with the following properties:
\begin{align}
 &x_i^\eps(t+2\pi\eps) = x_i^\eps(t)+1,&&
 y_i^\eps(t+2\pi\eps) = y_i^\eps(t)\label{eq:period}\\
 &x_1^\eps(t) = t/2\pi\eps+\phi_1^\eps(t/\eps),&&
 x_2^\eps(t) = t/2\pi\eps+\phi_2^\eps(t/\eps),\label{eq:form x}\\
 &y_1^\eps(t) = \eps\bigl(1+\eps\be_\eps+\psi_1^\eps(t/\eps)\bigr),&&
 y_2^\eps(t) = \eps\bigl(-1+\eps\be_\eps+\psi_2^\eps(t/\eps)\bigr).\label{eq:form y}
\end{align}
The correction functions $\phi_i^\eps,\psi_i^\eps$ satisfy
\begin{align}
 \|\phi_i^\eps\|_{H^1_{2\pi}}&\le C\eps^2,
 &\|\psi_i^\eps\|_{H^1_{2\pi}}&\le C\eps,\label{eq:bound1}\\
 \|\phi_1^\eps-\phi_2^\eps\|_{H^1_{2\pi}}&\le C\eps^3,
 &\|\psi_1^\eps-\psi_2^\eps\|_{H^1_{2\pi}}&\le C\eps^2.\label{eq:bound2}
\end{align}
The correction functions can be chosen to satisfy the normalizations
\begin{equation}\label{eq:normalization}
 \int_0^{2\pi}(\phi_1^\eps+\phi_2^\eps)\,\dd t=0,
 \qquad
 \int_0^{2\pi}(\psi_1^\eps+\psi_2^\eps)\,\dd t=0.
\end{equation}
Furthermore, $\be_\eps\to\be_0\in\R$ where $\be_0$ is explicitly defined in \eqref{eq:beta-limit-final}.
\end{Thm}

\begin{Rem}${}$\\ \indent
\rm 
a) Since $\Phi(x+1,y)=\Phi(x,y)$ property \eqref{eq:period} and the $2\pi$-periodicity of the correction functions imply $p^\eps(t+2\pi\eps)=p^\eps(t)$, hence $p^\eps$ has period $2\pi\eps$.

b) The operator $\cJ_\Ga$ depends only on the curvature of $\Si$ along $\Ga$, hence the nondegeneracy hypothesis depends only on this curvature function and the periodic boundary condition along $\Ga$, rather than on the geometry away from $\Ga$. This hypothesis is not satisfied on the round sphere due to its symmetries. On the round sphere however one can explicitly construct a family of solutions as in Theorem~\ref{thm:main}.

c) The constant $\be_0$ defined in \eqref{eq:beta-limit-final} depends not only on the curvature of $\Si$ along $\Ga$ but also on the Green function of the Laplace-Beltrami operator of $(\Si,g)$, hence it is an invariant depending on the global geometry, relative to the chosen normal field and the above normalizations.

d) The theorem holds for more general Hamiltonian systems. In fact, the Green function of the Laplace-Beltrami operator can be replaced by any function $G$ of the form
\[
  G(p_1,p_2) = -\frac{1}{2\pi}\log d_g(p_1,p_2) - H(p_1,p_2),
\]
where
\[
  H: \cN = \{(p_1,p_2)\in\Si\times\Si:d_g(p_1,p_2)<\rho\} \to \R
\]
is defined on a neighborhood of the diagonal, is smooth and symmetric: $H(p_1,p_2)=H(p_2,p_1)$. In particular, $\Si$ may have a boundary, provided that $\Ga$ is contained in its interior, and $G$ can be the Green Function of the Laplace-Beltrami operator with Dirichlet or Neumann boundary conditions. Also $\Si$ does not need to be compact. The Robin term must then be defined from $G$ with the same convention as in Section~\ref{sec:intro}.

e) Since the geodesic has length $1$ the solution $p^\eps(t)$ has minimal period $2\pi\eps$. Indeed, by \eqref{eq:period} the projection of either vortex trajectory onto $\Ga$ has winding number one over this period. If the orbit were an $m$-fold iterate of a shorter periodic orbit, its winding number would be divisible by $m$, hence $m=1$. If $\Ga$ has length $1/m$, the same construction, with arc-length coordinates and the nondegeneracy condition on $1/m$-periodic functions, yields a solution with minimal period $2\pi\eps/m$. In this case it is an interesting problem whether there also exists a solution with minimal period $2\pi\eps$, i.e.\ a subharmonic periodic dipole. Maybe Morse theoretic methods and an index theory for periodic dipoles would help so that one can distinguish between a periodic dipole with minimal period $2\pi\eps$ and an iterate of a periodic dipole with minimal period $2\pi\eps/m$. The index theory developed in \cite{Long:2002} for symplectic paths could be useful for this problem.

f) Theorem~\ref{thm:main} yields periodic dipoles near a prescribed nondegenerate closed geodesic, without imposing a symmetry of the surface. In addition to the question about the existence of subharmonic dipoles stated in e) one may ask many more interesting and challenging questions. For instance, can one place two or more dipoles on the surface that follow the same closed geodesic with a time shift? Such solutions exist on the round sphere. Can one place two or more periodic dipoles following distinct closed geodesics? We believe that the methods developed in this paper can be extended to address these questions.

g) A challenging problem is to regularize the periodic solutions from Theroem~\ref{thm:main} to obtain periodic solutions of the two-dimensional incompressible Euler equations or of the Gross-Pitaevskii equation on $\Si$. In a planar domain Davila, del Pino Musso and Wei \cite{DdPMW:2020} use gluing methods to construct solutions of the two-dimensional Euler flow with highly concentrated vorticities around moving vortices that follow on a finite time interval a solution of the point vortex dynamics. It is however an open problem to prove the existence of periodic solutions if the limit point vortex solution is periodic. We would also like to mention the construction of periodic vortex patch solutions of the Euler equations in symmetric settings, e.g.\ by Castro, C\'ordoba and G\'omez-Serrano \cite{Castro-Cordoba-Gomez:2019} in the plane, Hmidi, L.~Xue and Z.~Xue \cite{Hmidi-Xue-Xue:2026} in radial planar domains, Cao, Li and Wang \cite{Cao-Li-Wang:2025} on the round sphere. In \cite{Gelantalis-Sternberg:2012} Gelantalis and Sternberg  construct periodic solutions of the Gross-Pitaevskii equation
\[
  iU_t = \De U+\frac{1}{\eps^2}(1-|U|^2)U
\] 
on the unit sphere that possess $N$ vortices of degree $+1$ on a line of latitude and $N$ vortices of degree $-1$ placed symmetrically with respect to the equator. In the limit $\eps\to0$ these converge to periodic solutions of the $2N$-point vortex problem for $N$ dipoles.
\end{Rem}

\noindent\textit{Strategy of the proof of Theorem~\ref{thm:main}.}
We briefly describe the main steps of the argument. We first write the system in Fermi coordinates, obtaining a Hamiltonian system for $(z_1,z_2):\R \to \cS\times\cS$. After introducing the rescaled unknowns $u_i(t)=x_i(\eps t)$ and $v_i(t)=\eps^{-1}y_i(\eps t)$, the unperturbed dipole is given by two vortices travelling with the same tangential velocity, at normal positions $v_1=1$ and $v_2=-1$.  We therefore look for solutions of
the form
\[
 u_i(t) = \frac{t}{2\pi}+\phi_i(t),\qquad
 v_1(t) = 1+\eps\be+\psi_1(t),\qquad
 v_2(t) = -1+\eps\be+\psi_2(t),
\]
where the correction functions $\phi_i$, $\psi_i$ are $2\pi$-periodic, and $\be$ is a common normal displacement that is also to be determined.

A first essential ingredient is an expansion of $\log d_g^2$ in Fermi coordinates. In the regime relevant for the dipole, the normal
separation is of order $\eps$, whereas the tangential separation is much smaller, of order $\eps^3$; see \eqref{eq:form x}, \eqref{eq:form y}, and \eqref{eq:bound2}. Lemma~\ref{lem:distance-remainder} exploits this normal-dominated geometry to expand the singular part $\log d_g^2\big(\Phi(z_1),\Phi(z_2)\big)$ of the Hamiltonian and, crucially, of its first derivatives, with a remainder which is small enough also after differentiation with
respect to the scaled unknowns. This isolates the Euclidean dipole interaction and the first curvature correction. The regular part of the Hamiltonian is treated
separately.  Its special form, the symmetry with respect to permuting $z_1$ and $z_2$, and the dipole feature $\ka_1=-\ka_2$ of the two vortices imply the
diagonal identities of Lemma~\ref{lem:F-identities} which are responsible for several very helpful cancellations.

We then split the equations for $\phi_i$, $\psi_i$ and $\be$ into a principal linear part $L^\eps$ and a nonlinear
remainder. The kernel of $L^\eps$ is two-dimensional, consisting of the common translations of $\phi_1$, $\phi_2$, i.e.\ translations parallel to the geodesic, and the common translations of $\psi_1$, $\psi_2$, i.e.\ translations normal to the geodesic. We use an augmented reduction: the mean-value conditions \eqref{eq:normalization} fix the tangential phase and separate the constant normal displacement $\be$, while auxiliary scalar parameters account for the range conditions. The Hamiltonian structure makes one of the resulting compatibility conditions redundant.

Finally we pass to the sum--difference variables $\phi_1\pm\phi_2$ and $\psi_1\pm\psi_2$. After a suitable scaling of these we are able to obtain a limit equation for $\eps\to0$. A subtle and intricate contraction argument then produces a solution for all sufficiently small $\eps$, and the reconstruction in the original variables yields the one-parameter family of periodic dipole orbits asserted in Theorem~\ref{thm:main}.

\section{Preliminaries}\label{sec:prelim}

For $(p_1,p_2)\in\cF_2\Si$ with $d_g(p_1,p_2)<\rho$ small the Green function can be decomposed into a singular and a regular part:
\[
  G_\Si(p_1,p_2) = -\frac1{4\pi}\ln d_g^2(p_1,p_2) - H_\Si(p_1,p_2)
\]
where $H_\Si:\{(p_1,p_2)\in\Si\times\Si:d_g(p_1,p_2)<\rho\}\to\R$ is defined and smooth near the diagonal. Moreover $H_\Si$ is symmetric, and the Robin function satisfies $R_\Si(p)=H_\Si(p,p)$. Thus the Hamiltonian $\cH$ from \eqref{def:H} can be written as
\[
  \cH(p_1,p_2) = -\frac1{2\pi}\ln d_g^2(p_1,p_2) + H_\Si(p_1,p_1) + H_\Si(p_2,p_2) - 2H_\Si(p_1,p_2).
\]
In Fermi coordinates the metric $g$ has the form
\begin{equation}\label{metric}
  g(x,y)=\begin{pmatrix}  \la^2(x,y)&0 \\ 0&1 \end{pmatrix}
\end{equation}
where $\la$ satisfies
\begin{equation}\label{eq:lambda-synge}
  \la(x,y) = 1-\frac12\ka(x)y^2-\frac16\ka_\nu(x)y^3+O(y^4)\quad\text{as $y\to0$;}
\end{equation}
here $\ka(x)$ is the Gauss curvature of $\Si$ at $\ga(x)=\Phi(x,0)$ and $\ka_\nu(x)=\pa_{\nu(x)}K(\ga(x))$ is the derivative of $K$ at $\ga(x)$ in the normal direction. The expansion follows from
\[
 \pa_y^2\la(x,y)+K(\Phi(x,y))\la(x,y)=0,
 \qquad \la(x,0)=1,\quad \pa_y\la(x,0)=0.
\]
Indeed, $\pa_y^2\la(x,0)=-\ka(x)$ and $\pa_y^3\la(x,0)=-\ka_\nu(x)$. If $h:U\to\R$ is a differentiable function, 
then the gradient $\nabla h(\Phi(x,y))\in T_{\Phi(x,y)} \Si$ in Fermi coordinates is given via
\begin{equation}\label{eq:gradient}
  \nabla^g (h\circ \Phi)(x,y) = \left(\frac{1}{\la^2(x,y)}\pa_x (h\circ \Phi)(x,y), \pa_y (h\circ \Phi)(x,y)\right).
\end{equation}
The area form $\Om$ on $(\Si,g)$ induces an almost complex structure $J:T\Si\to T\Si$ defined by $\Om(v,Jw) = \langle v,w\rangle_g$. It is immediate to check that $J$ is a symplectic map and $J^{-1} = -J$. In Fermi coordinates the area form is given by
\[
  \Om_{(x,y)}(v,w) = \det\begin{pmatrix} \la(x,y) v_1& \la(x,y) w_1 \\ v_2&  w_2 \end{pmatrix} = \la(x,y)(v_1w_2-v_2w_1)
\]
and the almost complex structure by
$$ J_{(x,y)} = \begin{pmatrix}  0 & -\displaystyle{\frac{1}{\la(x,y) }} \\ \la(x,y) & 0\end{pmatrix}$$
i.e.\ 
\begin{equation}\label{eq:ac structure}
  J_{(x,y)}(v,w) = \left(-\frac{w}{\la(x,y)},\,\la(x,y)v\right).
\end{equation}
The almost complex structure on $T(\Si\times\Si)$ induced by $\om = \pi_1^*(\Om)+\pi_2^*(\Om)$ is in Fermi coordinates given by
\[
  J_{(z_1,z_2)}(v_1,w_1,v_2,w_2) = \bigl(J_{z_1}(v_1,w_1),\,J_{z_2}(v_2,w_2)\bigr).
\]

Now we introduce the functions
\begin{equation}\label{eq:synge}
  D_g:\cS\times\cS \to \R, \quad D_g(z_1,z_2) := d_g^2\big(\Phi(z_1),\Phi(z_2)\big),
\end{equation}
and
\begin{equation}\label{eq:F}
  F(z_1,z_2) := H_\Si\big(\Phi(z_1),\Phi(z_1)\big) + H_\Si\big(\Phi(z_2),\Phi(z_2)\big) - 2H_\Si\big(\Phi(z_1),\Phi(z_2)\big).
\end{equation}
The function $F$ is defined on the neighborhood
\[
  \mathcal V_\rho:=\{(z_1,z_2)\in\cS^2: d_g(\Phi(z_1),\Phi(z_2))<\rho\}
\]
of the lifted diagonal. For the Hamiltonian system we restrict to
\[
  \mathcal U_\rho:=\{(z_1,z_2)\in\cS^2: 0<d_g(\Phi(z_1),\Phi(z_2))<\rho\}.
\]
Thus collisions are excluded on the surface itself, including those represented by distinct lifts in the periodic coordinate $x$.

The system \eqref{HS} can be written in Fermi coordinates with the Hamiltonian
\begin{equation}\label{eq:def-HPhi} 
  \cH_\Phi:\mathcal U_\rho\to\R,\quad \cH_\Phi(z_1,z_2) := \cH\big(\Phi(z_1),\Phi(z_2)\big) = -\frac{1}{2\pi}\log D_g(z_1,z_2)+F(z_1,z_2)
\end{equation}
as
\begin{equation}\label{HS-Fermi}
\left\{
\begin{aligned}
  \dot z_1 &= J_{z_1} \nabla^g_{z_1}\cH_\Phi(z_1,z_2)\\
  \dot z_2 &= -J_{z_2} \nabla^g_{z_2}\cH_\Phi(z_1,z_2)
\end{aligned}
\right.
\end{equation} 
Then using \eqref{eq:gradient}--\eqref{eq:ac structure}, for $p_1,p_2$ near $\Ga$ the Hamiltonian system \eqref{HS} reduces to the system
\begin{equation}\label{eq:original-system}
\left\{
\begin{aligned}
 \dot x_1&=-\frac1{\la(x_1,y_1)}\partial_{y_1}\left(-\frac1{2\pi}\log D_g+F\right)(x_1,y_1,x_2,y_2),\\
 \dot y_1&= \frac1{\la(x_1,y_1)}\partial_{x_1}\left(-\frac1{2\pi}\log D_g+F\right)(x_1,y_1,x_2,y_2),\\
 \dot x_2&= \frac1{\la(x_2,y_2)}\partial_{y_2}\left(-\frac1{2\pi}\log D_g+F\right)(x_1,y_1,x_2,y_2),\\
 \dot y_2&=-\frac1{\la(x_2,y_2)}\partial_{x_2}\left(-\frac1{2\pi}\log D_g+F\right)(x_1,y_1,x_2,y_2).
\end{aligned}
\right.
\end{equation}

In the next section we investigate $\log D_g(z_1,z_2)$ and its derivatives in the special regime appearing in Theorem~\ref{thm:main}. We conclude this section with some algebraic identities for $F$ that simplify its Taylor expansion.

\begin{Lem}\label{lem:F-identities}
For every $z\in\cS$ there holds:
\[
 F(z,z)=0,
 \qquad
 \nabla_{z_1}F(z,z)=\nabla_{z_2}F(z,z)=0.
\]
Moreover, setting $A(z)=D^2_{z_1z_1}F(z,z)$ then
\begin{equation}\label{eq:F-hessian-block}
 D^2F(z,z)=
 \begin{pmatrix}
 A(z)&-A(z)\\
 -A(z)&A(z)
 \end{pmatrix}.
\end{equation}
In particular the second derivatives of $F$ satisfy on the diagonal:
\begin{equation}\label{eq:F-second-identities}
 F_{x_1x_1}=F_{x_2x_2}=-F_{x_1x_2},\qquad
 F_{y_1y_1}=F_{y_2y_2}=-F_{y_1y_2},
\end{equation}
and
\begin{equation}\label{eq:F-mixed-identities}
 F_{x_1y_1}=F_{x_2y_2}=-F_{x_1y_2}=-F_{x_2y_1}.
\end{equation}
\end{Lem}

\begin{proof}
The first identity follows directly from \eqref{eq:F}.  The symmetry
$F(z_1,z_2)=F(z_2,z_1)$ and the identity $F(z,z)\equiv0$ imply the vanishing of the
first derivatives on the diagonal.  Differentiating $F(z,z)\equiv0$ twice and using
symmetry gives the block form \eqref{eq:F-hessian-block}.
\end{proof}

As a consequence we obtain, setting $m=(z_1+z_2)/2$:
\begin{equation}\label{eq:F-gradient-structure}
\left\{
\begin{aligned}
  \nabla_{z_1}F(z_1,z_2) &= A(m)(z_1-z_2)+O(|z_1-z_2|^2),\\
  \nabla_{z_2}F(z_1,z_2) &= -A(m)(z_1-z_2)+O(|z_1-z_2|^2),
\end{aligned}
\right.
\qquad\text{as $|z_1-z_2|\to0$.}
\end{equation}
In particular,
\begin{equation}\label{eq:F-gradient-sum}
  \nabla_{z_1}F(z_1,z_2)+\nabla_{z_2}F(z_1,z_2) = O(|z_1-z_2|^2)\qquad\text{as $|z_1-z_2|\to0$.}
\end{equation}
There also holds for $|z_1-z_2|\to0$, with $e_1=(1,0)\in\R^2$:
\begin{align}
 \pa_{y_1}F(z_1,z_2)+\pa_{y_2}F(z_1,z_2) &= \frac12\,\big\langle\pa_yA(m)(z_1-z_2),z_1-z_2\big\rangle + O(|z_1-z_2|^4),
 \label{eq:Fysum-centered}\\
 \pa_{x_1}F(z_1,z_2)-\pa_{x_2}F(z_1,z_2) &= 2\big\langle A(m)(z_1-z_2),e_1\big\rangle + O(|z_1-z_2|^3),
 \label{eq:Fxdiff-centered}\\
 \pa_{x_1}F(z_1,z_2)+\pa_{x_2}F(z_1,z_2) &=\frac12\big\langle\pa_xA(m)(z_1-z_2),z_1-z_2\big\rangle + O(|z_1-z_2|^4).
 \label{eq:Fxsum-centered}
\end{align}
In order to see these set $d=z_1-z_2$, and observe that Lemma~\ref{lem:F-identities} and Taylor's formula imply
\begin{equation}\label{eq:F-centered}
 F\left(m+\frac d2,m-\frac d2\right) = \frac12\langle A(m)d,d\rangle+O(|d|^4).
\end{equation}
To justify the differentiation of the remainder, define $\widehat F(m,d):=F(m+d/2,m-d/2)$. The symmetry of $F$ implies that
$\widehat F(m,d)$ is even in $d$. Taylor's formula with integral remainder
therefore gives
\[
  \widehat F(m,d)=\frac12\langle A(m)d,d\rangle+R(m,d), \quad
  R=O(|d|^4),\quad \nabla_m R=O(|d|^4),\quad \nabla_d R=O(|d|^3).
\]
These estimates are uniform for $m$ in a fixed smaller tubular neighborhood,
modulo the period in its tangential coordinate. Using
\[
 \nabla_{z_1}=\frac12\nabla_m+\nabla_d,
 \qquad
 \nabla_{z_2}=\frac12\nabla_m-\nabla_d,
\]
we obtain \eqref{eq:Fysum-centered}--\eqref{eq:Fxsum-centered}, including the stated differentiated remainder estimates.
These identities imply several important cancellations when performing the Lyapunov-Schmidt reduction.

We can now define the limit $\be_0$ appearing in Theorem~\ref{thm:main}. Recall that $\ka_\nu(x)=\pa_{\nu(x)}K(\ga(x))$ denotes the normal derivative of the Gauss curvature along $\Ga$. For $z\in\cS$ we set
\begin{equation}\label{def:beta_0}
  A(z) = D^2_{z_1z_1}F(z,z)
            = \begin{pmatrix} a(z) &b(z)\\ b(z) &c(z) \end{pmatrix}
\end{equation}
Since by assumption the kernel of the operator $\cJ_\Ga$ defined in \eqref{def:J_Ga} is trivial there exists a unique $1$-periodic solution $Z_0$ of
\begin{equation}\label{eq:Zsolution}
  -Z_0''(x) -\ka(x) Z_0(x) = \frac13\ka_\nu(x) - 4\pi\big(\pa_yc(x,0)+2\pa_xb(x,0)\big).
\end{equation}
Then we have
\begin{equation}\label{eq:beta-limit-final}
 \beta_0 = \frac12\int_0^1 Z_0(x)\,\dd x.
\end{equation}

\section{Expansion of the squared distance and its derivatives}\label{sec:Taylor}

We fix $R>0$ and define for $\eps>0$ the set
\begin{equation}\label{eq:def-cA}
\begin{aligned}
  \cA_{R,\eps} &:= \big\{(z_1,z_2)\in\cS\times\cS:&&|x_i-\xi|\le R\eps^2 \text{ for some $\xi\in\R$},\ |x_1-x_2|\le R\eps^3, \\
               &  &&|y_1-\eps|\le R\eps^2,\ |y_2+\eps|\le R\eps^2\big\}
\end{aligned}
\end{equation}
There exists $c>0$, independent of $\eps$ for $\eps>0$ small, such that for $(z_1,z_2)\in\cA_{R,\eps}$:
\begin{equation}\label{eq:eta-cone}
  \frac1c < \frac{y_1-y_2}{\eps} <c,\quad \frac1c <\frac{\max\{|y_1|,|y_2|\}}{\eps} < c,\quad\text{and}\quad\frac{|x_1-x_2|}{|y_1-y_2|} <c\eps^2.
\end{equation}

The function $D_g(z_1,z_2) = d_g^2\big(\Phi(z_1),\Phi(z_2)\big)$, defined in \eqref{eq:synge} for $z_i\in \cS=\R\times(-\de,\de)$, is smooth in $\cA_{R,\eps}$ for $\eps>0$ small. Below we will use the notation 
\begin{equation}\label{notation1}
  m = \frac{x_1+x_2}2,\quad X = x_1-x_2,\quad Y = y_1-y_2, \quad \eta = \max\{|y_1|,|y_2|\}.
\end{equation}

\begin{Lem}\label{lem:synge1}
  For $(z_1,z_2)\in\cA_{R,\eps}$ with $|x_i-\xi|\le R\eps^2$ there holds:
  \[
    D_g(z_1,z_2) = (x_1-x_2)^2+(y_1-y_2)^2 -\frac{\ka(\xi)}3(x_1-x_2)^2(y_1^2+y_1y_2+y_2^2) + \cR(m,y_1,y_2,X;\xi)
  \]
  Here $\xi$ is held fixed in all partial derivatives below. The remainder term $\cR=\cR(m,y_1,y_2,X;\xi)$ satisfies the following estimates on $\cA_{R,\eps}$ with a constant $C>0$ that depends on $R$ but is independent of $\xi$ and $\eps$ for $\eps>0$ small:
  \begin{align}
   &|\cR| \le C\big(\eta^3X^2+|X|^3\big) \label{eq:R1}\\
   &\big|\pa_X\cR\big| \le C\eta^3|X| &&\big|\pa_{y_i}\cR\big| \le C\eta^2X^2,\label{eq:R2}\\
   &\big|\pa_X^2\cR\big| \le C\eta^3 &&\big|\pa_X\pa_{y_i}\cR\big| \le C\eta^2|X| &&& \big|\pa_{y_i}\pa_{y_j}\cR\big| \le C\eta X^2,\label{eq:R3}\\
   &\big|\pa_m\cR\big| \le C\big(\eta^2X^2+|X|^3\big) &&\big|\pa_m\pa_X\cR\big| \le C\big(\eta^2|X|+X^2\big),\label{eq:R4}\\
   &\big|\pa_m\pa_{y_i}\cR\big| \le C\big(\eta X^2+|X|^3\big) &&\big|\pa_m^2\cR\big| \le C\big(\eta^2X^2+|X|^3\big).\label{eq:R5}
  \end{align}
\end{Lem}

\begin{proof}
For fixed $(m,y_1,y_2)$ define
\begin{equation}\label{eq:D-of-X}
  \cD(X) := D_g\bigl((m+X/2,y_1),(m-X/2,y_2)\bigr).
\end{equation}
For $X=0$, i.e.\ $x_1=x_2=m$, the unique minimizing geodesic is the normal segment $x\equiv m$, hence
\begin{equation}\label{eq:D-zero-first}
 \cD(0) = Y^2,
 \qquad
 \cD'(0)=0.
\end{equation}
We now compute $\cD''(0)$ explicitly in Fermi coordinates.

Recall that the squared distance between two points is the minimum of the energy
\[
  \cJ[x,y] = \int_0^1\bigl(\lambda(x,y)^2\dot x^2+\dot y^2\bigr)\,\dd\theta
\]
among curves $(x,y):[0,1]\to\cS$ joining the endpoints $z_1,z_2$, with the minimum attained by the affinely parametrized minimizing geodesic. At $X=0$ the minimizer is
\[
 \bar x(\theta)=m,
 \qquad
 \bar y(\theta)=y_2+Y\theta.
\]
Let $(x_X,y_X)$ be the affinely parametrized minimizing geodesic joining $(m-X/2,y_2)$ to $(m+X/2,y_1)$. 
For nearby endpoints it depends smoothly on $X$. Write
\[
 x_X=m+Xh+O(X^2),\qquad y_X=\bar y+Xk+O(X^2),
\]
where 
\[
 h(0)=-\tfrac12,\quad h(1)=\tfrac12,\qquad k(0)=k(1)=0,
\]
and the expansions hold in $C^1([0,1])$. Linearizing the geodesic equations along $(m,\bar y)$ gives
\[
 \frac{\dd}{\dd\theta}\bigl(\la^2(m,\bar y)\dot h\bigr)=0,
 \qquad \ddot k=0.
\]
This implies $k=0$. In the expansion of the energy, the term involving the second
variation of $y_X$ integrates to zero, since $\dot{\bar y}=Y$ is constant
and the endpoints of $y_X$ are independent of $X$. Consequently,
\[
 \cD(X)=\cJ[x_X,y_X]
 =Y^2+X^2\int_0^1\la^2(m,\bar y)\dot h^2\,\dd\theta+o(X^2).
\]
The function $h$ solves the Euler--Lagrange equation of the strictly convex
quadratic functional defined in \eqref{eq:B-variational} below and hence is its unique minimizer. It follows that
\begin{equation}\label{eq:D-second}
  \frac12\cD''(0) = B(m,y_1,y_2),
\end{equation}
where $B(m,y_1,y_2)$ is the minimum of
\begin{equation}\label{eq:B-variational}
 h\mapsto \int_0^1\la^2(m,\bar y)\dot h^2\,\dd\theta
\end{equation}
under the above endpoint conditions for $h:[0,1]\to\R$. The Euler--Lagrange equation for the minimum is
\begin{equation}\label{eq:jacobi-fermi-h}
 \frac{\dd}{\dd\theta} \left(\la^2(m,y_2+Y\theta)\dot h(\theta)\right)=0.
\end{equation}
This is the scalar Jacobi equation for the tangential variation of the geodesic $(\bar x,\bar y)$.  Now from \eqref{eq:jacobi-fermi-h} we deduce
\[
 \dot h(\theta) = C\la^{-2}(m,y_2+Y\theta).
\]
Therefore $h(1)-h(0)=1$ implies
\[
  C=\left( \int_0^1\la^{-2}(m,y_2+Y\theta)\,\dd\theta \right)^{-1}.
\]
Substitution into \eqref{eq:B-variational} gives the exact formula
\begin{equation}\label{eq:B-exact}
 B(m,y_1,y_2)
 =\left( \int_0^1\la^{-2}(m,y_2+Y\theta)\,\dd\theta \right)^{-1}
 =\frac{Y}{\displaystyle \int_{y_2}^{y_1}\la^{-2}(m,y)\,\dd y}.
\end{equation}
The last quotient extends smoothly across $Y=0$.

Next we expand this quotient. From \eqref{eq:lambda-synge} we deduce
\begin{align*}
 \frac1Y\int_{y_2}^{y_1}\la(m,y)^{-2}\,\dd y
 &=1+\frac{\ka(m)}{3} \big(y_1^2+y_1y_2+y_2^2\big)\\
 &\hspace{.75cm}
     +\frac{\ka_\nu(m)}{12}(y_1^3+y_1^2y_2+y_1y_2^2+y_2^3) + O(\eta^4).
\end{align*}
Taking the reciprocal gives
\begin{equation}\label{eq:B-expansion}
\begin{aligned} B(m,y_1,y_2)
 &=1-\frac{\ka(m)}3
 \big(y_1^2+y_1y_2+y_2^2\big)\\
 &\hspace{.75cm}
 -\frac{\ka_\nu(m)}{12}
 \big(y_1^3+y_1^2y_2+y_1y_2^2+y_2^3\big) + O(\eta^4).
\end{aligned}
\end{equation}
The remainder $E_4$ denoted by $O(\eta^4)$ in \eqref{eq:B-expansion} is a smooth Taylor remainder. For $a+b\le2$ and $0\le j\le2$ it satisfies
$|\pa_m^j\pa_{y_1}^a\pa_{y_2}^b E_4|\le C\eta^{4-a-b}$,
uniformly in $m$ modulo its period. Since $\cD$ is smooth, Taylor's formula in the single variable $X$ gives, uniformly 
for $(m,y_1,y_2)$ in $\cA_{R,\eps}$:
\begin{equation}\label{eq:D-Taylor-X}
 D_g(z_1,z_2) = Y^2 + B(m,y_1,y_2)X^2 + X^3R_3(m,y_1,y_2,X),
\end{equation}
where $R_3$ and all derivatives of $R_3$ are uniformly bounded. 

Replacing $\ka(m)$ by $\ka(\xi)$ in \eqref{eq:B-expansion}, using $|m-\xi|=O(\eps^2)=O(\eta^2)$, and absorbing the cubic normal term and the $X^3R_3$ term into the remainder, we obtain
\begin{equation}\label{eq:distance}
  D_g(z_1,z_2) = X^2 + Y^2 - \frac{\ka(\xi)}3X^2\big(y_1^2+y_1y_2+y_2^2\big) + \cR(m,y_1,y_2,X;\xi),
\end{equation}
where $\cR$ satisfies \eqref{eq:R1}. Differentiating \eqref{eq:D-Taylor-X} and \eqref{eq:B-expansion} and using $|X|\le C\eta^3$ gives the estimates \eqref{eq:R2}--\eqref{eq:R5}.
\end{proof}

We now pass to the expansion of $\log D_g$, using the notation in \eqref{notation1}.

\begin{Lem}\label{lem:distance-remainder}
  For the fixed $R>0$, there exist $\eps_0>0$ and $C>0$ such that for $0<\eps<\eps_0$ and $(z_1,z_2)\in\cA_{R,\eps}$ with $|x_i-\xi|\le R\eps^2$ there holds:
  \begin{equation}\label{eq:log-distance}
    \log D_g(x_1,y_1,x_2,y_2) = \log(X^2+Y^2) - \frac{\ka(\xi)}3\cdot\frac{X^2(y_1^2+y_1y_2+y_2^2)}{X^2+Y^2} + \cR_{log}
  \end{equation}
  where the remainder term $\cR_{log}$ satisfies:
  \begin{align}
    &\cR_{\log} = \frac{1}{X^2+Y^2}\cdot\left(\cR + \left(\frac{\ka(\xi)}3\cdot X^2(y_1^2+y_1y_2+y_2^2)-\cR\right)^2\cdot O(1)\right)
           &&\text{as $\eps\to0$} \label{eq:Rlog-structure}\\
    &\big|\pa_{x_i}\cR_{log}\big| \le C\eps^4, &&\hspace{-8cm}\big|\pa_{y_i}\cR_{log}\big| \le C\eps^6 \label{eq:distance-remainder-strong}\\
    &\big|\pa_{x_i x_j}^2\cR_{\log}\big|\le C\eps, &&\hspace{-8cm}\big|\pa_{x_i y_j}^2\cR_{\log}\big| \le C\eps^3,
          \qquad\big|\pa_{y_i y_j}^2\cR_{\log}\big| \le C\eps^5 \label{eq:Rlog-second-derivatives}
  \end{align}
  The remainder also depends on $\xi$, which is held fixed in every partial derivative. The estimates are uniform in $\xi\in\R$.
\end{Lem}

\begin{proof}
 In order to shorten notation we consider the smooth function 
\[
 h:(-1,\infty)\to\R,\qquad
 h(u)=\begin{cases}\displaystyle\frac{\log(1+u)-u}{u^2},&u\ne0,\\[4pt]
 -\tfrac12,&u=0.\end{cases}
\]
and set
\[
  Q=X^2+Y^2,\quad S=y_1^2+y_1y_2+y_2^2,\quad M_4=\frac{\ka(\xi)}3\cdot X^2S,\quad u=\frac{-M_4+\cR}{Q}.
\]
Then we have the identities
\[
  D_g(z_1,z_2) = Q-M_4+\cR = Q\cdot\left(1-\frac{M_4}Q+\frac{\cR}Q\right).
\]
and
\[
 \log D_g=\log Q-\frac{M_4}{Q}+\frac{\cR}{Q}+u^2h(u).
\]
so that
\[
  \cR_{log} = \frac{\cR}{Q}+u^2h(u).
\]
Now observe that $u\to0$ as $\eps\to0$, hence \eqref{eq:log-distance} and \eqref{eq:Rlog-structure} follow from $h(0)\in\R$. 
Next using \eqref{eq:eta-cone} direct differentiation yields
\begin{equation}\label{eq:Rlog-first-eta}
 \big|\pa_X\cR_{\log}\big| \le C\eta|X| \quad\text{and}\quad \big|\pa_{y_i}\cR_{\log}\big| \le CX^2.
\end{equation}
In order to pass from $X=x_1-x_2$ to $x_1,x_2$, recall that
\[
  \pa_{x_1}=\frac12\pa_m+\pa_X \quad\text{and}\quad \pa_{x_2}=\frac12\pa_m-\pa_X.
\]
Now \eqref{eq:R4} implies
\[
 \big|\pa_{x_i}\cR_{\log}\big|\le C\eta|X| \quad\text{and}\quad \big|\pa_{y_i}\cR_{\log}\big| \le CX^2.
\]
Using \eqref{eq:eta-cone} again we obtain \eqref{eq:distance-remainder-strong}.

Concerning the second derivatives, the leading part of $\cR_{log}$ is $\cR/Q$. The estimates \eqref{eq:R2}, \eqref{eq:R3} yield
\[
 \big|\pa_X^2(\cR/Q)\big|\le C\eta,
 \qquad
 \big|\pa_X\pa_{y_i}(\cR/Q)\big|\le C|X|,
 \qquad
 \big|\pa_{y_i}\pa_{y_j}(\cR/Q)\big| \le C\frac{X^2}{\eta}.
\]
For instance, straightforward calculation gives
\[
  \pa_X^2\left(\frac{\cR}{Q}\right) = \frac{\pa_X^2\cR}{Q} - 2\frac{(\pa_X\cR)(\pa_XQ)}{Q^2} - \frac{\cR\,\pa_X^2Q}{Q^2} + 2\frac{\cR(\pa_XQ)^2}{Q^3},
\]
and each term is of order $O(\eta)$ because $\pa_X^2\cR = O(\eta^3)$, $\pa_X\cR=O(\eta^3|X|)$, $\cR=O(\eta^3X^2+|X|^3)$, $|X|\le C\eta^3$, and $Q/\eta^2$ is bounded away from $0$ and $\infty$ on $\cA_{R,\eps}$ uniformly in $\eps$ by \eqref{eq:eta-cone}. The mixed and $y$-derivative estimates of $\cR/Q$ are obtained in the same way.

Next observe that
\[
 \big|\pa_m\pa_X(\cR/Q)\big|\le C|X|,
 \qquad
 \big|\pa_m\pa_{y_i}(\cR/Q)\big|\le C\frac{X^2}{\eta},
 \qquad
 \big|\pa_m^2(\cR/Q)\big|\le CX^2.
\]
as a consequence of \eqref{eq:R4} and \eqref{eq:R5}. Finally the bounds above and direct differentiation give
\[
\begin{gathered}
 |u|\le CX^2,\quad |\pa_{x_i}u|\le C|X|,\quad
 |\pa_{y_i}u|\le CX^2/\eta,\\
 |\pa_{x_i x_j}^2u|\le C,\quad
 |\pa_{x_i y_j}^2u|\le C|X|/\eta,\quad
 |\pa_{y_i y_j}^2u|\le CX^2/\eta^2.
\end{gathered}
\]
Since $h$ and its first two derivatives are bounded, the chain rule yields for $v:=u^2h(u)$:
\[
  |\pa_{x_i}v|\le C|X|^3,\quad |\pa_{y_i}v|\le CX^4/\eta,
\]
and
\[
 |\pa_{x_i x_j}^2v|\le CX^2,\quad
 |\pa_{x_i y_j}^2v|\le C|X|^3/\eta,\quad
 |\pa_{y_i y_j}^2v|\le CX^4/\eta^2.
\]
These terms satisfy the required bounds because $|X|\le C\eta^3$.
Combining them with the estimates for $\cR/Q$ gives
\[
 \big|\pa_{x_i x_j}^2\cR_{\log}\big| \le C\eta,
 \qquad
 \big|\pa_{x_i y_j}^2\cR_{\log}\big| \le C|X|,
 \qquad
 \big|\pa_{y_i y_j}^2\cR_{\log}\big|\le C\frac{X^2}{\eta},
\]
which in turn implies \eqref{eq:Rlog-second-derivatives} in the region $\cA_{R,\eps}$.
\end{proof}

\begin{Rem}\label{rem:aux-distance-file}
The derivation uses the linearized geodesic equations along the normal segment to obtain the exact coefficient $B$. Its Taylor expansion yields the fourth-order curvature term, while Taylor expansion in the tangential separation $X$ and the bound $|X|\le C\eta^3$ give the sharper remainder estimates. 
\end{Rem}

\begin{Rem}
For the unit sphere, take the equator as the reference geodesic and use Fermi coordinates
\[
 \Phi(x,y)=(\cos y\cos x,\cos y\sin x,\sin y).
\]
Then
\[
 g=\cos^2y\,\dd x^2+\dd y^2, \qquad \ka\equiv1,
\]
and for $X=x_1-x_2$,
\begin{equation}\label{eq:sphere-cos}
 \cos d_g(z_1,z_2) = \sin y_1\sin y_2+\cos y_1\cos y_2\cos X.
\end{equation}
Expanding \eqref{eq:sphere-cos} at $(X,y_1,y_2)=(0,0,0)$ and using $\cos d=1-d^2/2+d^4/24+O(d^6)$ gives
\begin{equation}\label{eq:sphere-distance}
 d_g^2(z_1,z_2) = X^2+(y_1-y_2)^2 - \frac13X^2(y_1^2+y_1y_2+y_2^2) + O\bigl((|X|+|y_1|+|y_2|)^5\bigr).
\end{equation}
This verifies the coefficient in \eqref{eq:distance} for curvature $\ka=1$. For a sphere of constant curvature $\ka$, scaling gives the coefficient
$-\ka/3$. 
Similarly one can explicitly check the expansion on the hyperbolic plane. Knowing the expansion of $D_g$ in Lemma~\ref{lem:synge1} for surfaces of constant curvature one can derive the expansion in the general case by a comparison argument. We thank Alexander Lytchak for this argument. The comparison argument however does not yield the estimates \eqref{eq:R2}--\eqref{eq:R5} and \eqref{eq:Rlog-structure}--\eqref{eq:Rlog-second-derivatives} which are essential for the proof of Theorem~\ref{thm:main}.
\end{Rem}

\section{Proof of Theorem~\ref{thm:main}}\label{sec:proof}
\subsection{Ansatz and exact equations for the remainders}
In order to find the $2\pi\eps$-periodic solutions of \eqref{HS}, we consider system \eqref{eq:original-system} and rescale it by setting
\begin{equation}\label{eq:scaling}
 u_i(t)=x_i(\eps t) \quad\text{and}\quad v_i(t)=\frac1\eps y_i(\eps t).
\end{equation}
This leads to the scaled system
\begin{equation}\label{eq:scaled-system}
\left\{
\begin{aligned}
 \dot u_1&=\ga_1\left( \frac{\eps}{2\pi}\pa_{y_1}\log D_g-\eps\pa_{y_1}F\right)(u_1,\eps v_1,u_2,\eps v_2),\\
 \dot u_2&=\ga_2\left( -\frac{\eps}{2\pi}\pa_{y_2}\log D_g+\eps\pa_{y_2}F\right)(u_1,\eps v_1,u_2,\eps v_2),\\
 \dot v_1&=\ga_1\left( -\frac1{2\pi}\pa_{x_1}\log D_g+\pa_{x_1}F\right)(u_1,\eps v_1,u_2,\eps v_2),\\
 \dot v_2&=\ga_2\left( \frac1{2\pi}\pa_{x_2}\log D_g-\pa_{x_2}F\right)(u_1,\eps v_1,u_2,\eps v_2),
\end{aligned}
\right.
\end{equation}
where $F$ is defined in \eqref{eq:F}, and
\begin{equation}\label{eq:gamma-def}
 \ga_i=\frac1{\la(u_i,\eps v_i)}, \qquad i=1,2.
\end{equation}
The periodicity conditions are
\begin{equation}\label{eq:periodicity}
 u_i(t+2\pi)=u_i(t)+1 \quad\text{and}\quad v_i(t+2\pi)=v_i(t).
\end{equation}
We seek solutions of \eqref{eq:scaled-system}, \eqref{eq:periodicity} for $\eps>0$ small by setting $\ubar(t)=t/2\pi$ and making the ansatz
\begin{equation}\label{eq:ansatz}
\left\{
\begin{aligned}
 u_1&=\ubar+\phi_1,
 &\qquad
 u_2&=\ubar+\phi_2,\\
 v_1&=1+\beta\eps+\psi_1,
 &
 v_2&=-1+\beta\eps+\psi_2.
\end{aligned}
\right.
\end{equation}
The following computations hold provided
\begin{equation}\label{eq:phi-psi}
  |\phi_i|\le C\eps^2, \quad |\psi_i|\le C\eps,\quad |\phi_1-\phi_2|\le C\eps^3, \quad |\be|\le B
\end{equation}
for some fixed $B,C>0$. Then the points $(u_i(t),\eps v_i(t))$ lie in $\cA_{R,\eps}$ for $R=B+C$, with $\xi=\ubar(t)$. All remainder estimates below are uniform under these bounds. We have by \eqref{eq:lambda-synge}:
\begin{equation}\label{eq:gamma}
\left\{
\begin{aligned}
 \ga_1
 &=1+\frac12\ka(\ubar)\eps^2 + \ka(\ubar)\eps^2(\be\eps+\psi_1) + \frac16\ka_\nu(\ubar)\eps^3 + O(\eps^4),\\
 \ga_2
 &=1+\frac12\ka(\ubar)\eps^2 - \ka(\ubar)\eps^2(\be\eps+\psi_2) - \frac16\ka_\nu(\ubar)\eps^3 + O(\eps^4),
\end{aligned}
\right.
\end{equation}
so that
\begin{equation}\label{eq:gamma-difference}
 \ga_1-\ga_2 =\ka(\ubar)\eps^2\bigl(2\be\eps+\psi_1+\psi_2\bigr) + \frac13\ka_\nu(\ubar)\eps^3+O(\eps^4).
\end{equation}
Using the computation of $\log D_g$ in \eqref{eq:log-distance}, a straightforward computation shows that the system \eqref{eq:scaled-system} for $u_i$ and $v_i$ is equivalent to the following system for $\phi=(\phi_1,\phi_2)$ and $\psi=(\psi_1,\psi_2)$:
\begin{equation}\label{eq:remainder-system}
\left\{
\begin{aligned}
 L^\eps_{\phi1}(\phi,\psi) &:= \dot\phi_1 + \frac{\psi_1-\psi_2}{4\pi}\\
     &= A_{\phi} + B_{\phi1} + \eps\ga_1\cR_{\phi1} - \eps\ga_1\pa_{y_1}F(u_1,\eps v_1,u_2,\eps v_2) =: \cN^\eps_{\phi1}(\phi,\psi,\be)\\
 L^\eps_{\phi2}(\phi,\psi) &:= \dot\phi_2 + \frac{\psi_1-\psi_2}{4\pi}\\
     &= A_{\phi} + B_{\phi2} + \eps\ga_2\cR_{\phi2} + \eps\ga_2\pa_{y_2}F(u_1,\eps v_1,u_2,\eps v_2) =: \cN^\eps_{\phi2}(\phi,\psi,\be)\\
 L^\eps_{\psi1}(\phi,\psi) &:= \dot\psi_1 + \frac{\phi_1-\phi_2}{4\pi\eps^2} \\
     &= A_{\psi} + B_{\psi1} + \ga_1\cR_{\psi1} + \ga_1\pa_{x_1}F(u_1,\eps v_1,u_2,\eps v_2) =: \cN^\eps_{\psi1}(\phi,\psi,\be)\\
 L^\eps_{\psi2}(\phi,\psi) &:= \dot\psi_2 + \frac{\phi_1-\phi_2}{4\pi\eps^2} \\
     &= A_{\psi} + B_{\psi2} + \ga_2\cR_{\psi2} - \ga_2\pa_{x_2}F(u_1,\eps v_1,u_2,\eps v_2) =: \cN^\eps_{\psi2}(\phi,\psi,\be)
\end{aligned}
\right.
\end{equation}
Here $A_{\phi}$, $A_{\psi}$, $B_{\phi}$, $B_{\psi}$ are defined as follows:
\begin{equation}\label{eq:AB}
\left\{
\begin{aligned}
  &A_{\phi} = \frac{\eps^2}{\pi}\cdot\frac{2+\psi_1-\psi_2}{(\phi_1-\phi_2)^2+\eps^2(2+\psi_1-\psi_2)^2}-\frac1{2\pi}+\frac{\psi_1-\psi_2}{4\pi}\\
  &B_{\phi i} = (\ga_i-1)\cdot\frac{\eps^2}{\pi}\cdot\frac{2+\psi_1-\psi_2}{(\phi_1-\phi_2)^2+\eps^2(2+\psi_1-\psi_2)^2}\\
  &A_{\psi} = \frac1\pi\cdot\frac{\phi_2-\phi_1}{(\phi_1-\phi_2)^2+\eps^2(2+\psi_1-\psi_2)^2}+\frac{\phi_1-\phi_2}{4\pi\eps^2}\\
  &B_{\psi i} = (\ga_i-1)\cdot\frac1\pi\cdot\frac{\phi_2-\phi_1}{(\phi_1-\phi_2)^2+\eps^2(2+\psi_1-\psi_2)^2}
\end{aligned}
\right.
\end{equation}
Using the ansatz \eqref{eq:ansatz} and setting
\[
  Q = Q(\phi,\psi,\be) = (u_1-u_2)^2+\eps^2(v_1-v_2)^2,\quad S = S(\phi,\psi,\be) = \eps^2\big(v_1^2+v_1v_2+v_2^2\big)
\]
the remainder terms $\cR_{\phi i}= \cR_{\phi i}(\phi,\psi,\be)$ and $\cR_{\psi i} = \cR_{\psi i}(\phi,\psi,\be)$ are defined by:
\begin{equation}\label{eq:distance-corrections}
\left\{
\begin{aligned}
  \cR_{\phi1}
    &= \frac1{2\pi}\bigg[-\frac{\ka(\ubar)}3\cdot (u_1-u_2)^2\cdot\eps\cdot\bigg(\frac{2v_1+v_2}{Q}-\frac{2S(v_1-v_2)}{Q^2}\bigg)\\
    &\hspace{1cm} +\pa_{y_1}\cR_{log}(u_1,\eps v_1,u_2,\eps v_2)\bigg],\\
  \cR_{\phi2}
    &= -\frac1{2\pi}\bigg[-\frac{\ka(\ubar)}3\cdot(u_1-u_2)^2\cdot\eps\cdot\bigg(\frac{v_1+2v_2}{Q}+\frac{2S(v_1-v_2)}{Q^2}\bigg)\\
    &\hspace{1cm} +\pa_{y_2}\cR_{log}(u_1,\eps v_1,u_2,\eps v_2)\bigg],\\
  \cR_{\psi1}
    &= -\frac1{2\pi}\bigg[-\frac{2\ka(\ubar)}{3}\cdot\eps^2\cdot\frac{(u_1-u_2)S(v_1-v_2)^2}{Q^2}+\pa_{x_1}\cR_{log}(u_1,\eps v_1,u_2,\eps v_2)\bigg],\\
  \cR_{\psi2}
    &= \frac1{2\pi}\bigg[\frac{2\ka(\ubar)}{3}\cdot\eps^2\cdot\frac{(u_1-u_2)S(v_1-v_2)^2}{Q^2}+\pa_{x_2}\cR_{log}(u_1,\eps v_1,u_2,\eps v_2)\bigg].
\end{aligned}
\right.
\end{equation}
In these formulas, the derivatives of $\cR_{\log}$ are taken with the parameter $\xi$ fixed, and then evaluated at $\xi=\ubar(t)$. The dependence on $(\phi,\psi,\be)$ is through the ansatz for all four variables $u_i,v_i$.
The right hand side of \eqref{eq:remainder-system} still contains the contributions from the derivatives of $F$, including their linear terms. We leave these contributions in the right hand side rather than incorporating them into $L^\eps$. By Lemma~\ref{lem:F-identities}, the leading linearization of $\nabla F$ on the diagonal acts on the difference
\[
 z_1-z_2=\bigl(\phi_1-\phi_2,\,2\eps+\eps(\psi_1-\psi_2)\bigr),
 \qquad z_i=(u_i,\eps v_i).
\]
Thus it includes forcing terms as well as terms involving the differences of the remainders. Its coefficients depend on the midpoint $(z_1+z_2)/2$, so the full contributions need not depend only on those differences. In the sequel we also impose the bound $\psi_1-\psi_2=O(\eps^2)$ in addition to the bounds in \eqref{eq:phi-psi}.

\subsection{The principal linear operator}
Let $H^1_{2\pi}$ and $L^2_{2\pi}$ be the spaces of $2\pi$-periodic $H^1_{loc}$-, respectively $L^2_{loc}$-functions $\R\to\R$. It is easy to see that the kernel of the operator
\begin{equation}\label{eq:L}
  L^\eps = \big(L^\eps_{\phi1},L^\eps_{\phi2},L^\eps_{\psi1},L^\eps_{\psi2}\big):(H^1_{2\pi})^4\to(L^2_{2\pi})^4
\end{equation}
is given by
\begin{equation}\label{eq:kernel}
 K = \ker(L^\eps) = \operatorname{span}\{(1,1,0,0),(0,0,1,1)\}.
\end{equation}
Clearly we have
\begin{equation}\label{eq:Kperp}
  K^\perp = \left\{(\phi_1,\phi_2,\psi_1,\psi_2)\in (H^1_{2\pi})^4:\int_0^{2\pi}(\phi_1+\phi_2)\,\dd t=0,\  \int_0^{2\pi}(\psi_1+\psi_2)\,\dd t=0.\right\}
\end{equation}
The range of $L^\eps$ can also be explicitly computed.
\begin{Prop}\label{prop:linear}
For $(f_1,f_2,g_1,g_2)\in(L^2_{2\pi})^4$ the equation
\[
 L^\eps(\phi_1,\phi_2,\psi_1,\psi_2) = (f_1,f_2,g_1,g_2)
\]
has a solution in $(H^1_{2\pi})^4$ if and only if
\begin{equation}\label{eq:range}
 \int_0^{2\pi}(f_1-f_2)\,\dd t=0 \quad\text{and}\quad \int_0^{2\pi}(g_1-g_2)\,\dd t=0.
\end{equation}
There is a unique solution in $K^\perp$. The following estimates hold for the solution with a constant $C$ independent of $\eps>0$:
\begin{align}
 \|\phi_1-\phi_2\|_{H^1}
 &\le C\left(\|f_1-f_2\|_{L^2}
 +\eps^2\|g_1+g_2\|_{L^2}\right),\label{eq:lin-est-p}\\
 \|\psi_1-\psi_2\|_{H^1}
 &\le C\left(\|g_1-g_2\|_{L^2}
 +\|f_1+f_2\|_{L^2}\right),\label{eq:lin-est-q}\\
 \|\phi_1+\phi_2\|_{H^1}
 &\le C\left(\|f_1+f_2\|_{L^2}+\|g_1-g_2\|_{L^2}\right),\label{eq:lin-est-phisum}\\
 \|\psi_1+\psi_2\|_{H^1}
 &\le C\left(\|g_1+g_2\|_{L^2}
 +\eps^{-2}\|f_1-f_2\|_{L^2}\right).\label{eq:lin-est-psisum}
\end{align}
\end{Prop}

\begin{proof}
Subtracting the first two equations and the last two equations gives
\[
 \dot\phi_1-\dot\phi_2 = f_1-f_2 \quad\text{and}\quad \dot\psi_1-\dot\psi_2 = g_1-g_2,
\]
which yields \eqref{eq:range}. Integrating the sums of the first two and the last two equations gives the mean values
\begin{equation}\label{eq:means-pq}
 \int_0^{2\pi}(\phi_1-\phi_2)\,\dd t = 2\pi\eps^2\int_0^{2\pi}(g_1+g_2)\,\dd t \quad\text{and}\quad
 \int_0^{2\pi}(\psi_1-\psi_2)\,\dd t = 2\pi\int_0^{2\pi}(f_1+f_2)\,\dd t.
\end{equation}
In order to prove sufficiency, assume \eqref{eq:range} and set
\[
 p=\phi_1-\phi_2,\quad q=\psi_1-\psi_2,\qquad
 a=\phi_1+\phi_2,\quad b=\psi_1+\psi_2.
\]
The system is equivalent to
\[
\begin{aligned}
 \dot p&=f_1-f_2,&\qquad \dot q&=g_1-g_2,\\
 \dot a&=f_1+f_2-\frac{q}{2\pi},&
 \dot b&=g_1+g_2-\frac{p}{2\pi\eps^2}.
\end{aligned}
\]
The first two equations have solutions $p,q\in H^1_{2\pi}$ which are uniquely determined by prescribing the integrals in \eqref{eq:means-pq}. Then  the right hand sides of the equations for $a,b$ have zero mean. They therefore admit unique primitives $a,b\in H^1_{2\pi}$ with zero mean. Recovering
\[
  \phi_1=(a+p)/2,\quad \phi_2=(a-p)/2,\qquad
  \psi_1=(b+q)/2,\quad \psi_2=(b-q)/2
\]
gives the unique solution in $K^\perp$. The Poincar\'e inequality applied to $p,q$ after subtracting their means, and then to the zero-mean functions $a,b$, yields \eqref{eq:lin-est-p}--\eqref{eq:lin-est-psisum}, with constants independent of $\eps$.
\end{proof}

\begin{Rem}
There is no $\eps$-uniform bound for $(L^\eps)^{-1}$ in the standard unweighted $(H^1_{2\pi})^4$ norm for arbitrary data because the term $\eps^{-2}(f_1-f_2)$ in \eqref{eq:lin-est-psisum} is unavoidable. The refined estimates above are exactly what is needed in the sequel.
\end{Rem}

\subsection{Bounds for the nonlinear right-hand side}
We want to solve the remainder equation \eqref{eq:remainder-system} $L^\eps(\phi,\psi)=\cN^\eps(\phi,\psi,\be)$ in the domain
\begin{equation}\label{eq:CR}
\begin{aligned}
\cK_{M,\eps}&=\big\{(\phi,\psi): |\phi_1|+|\phi_2|\le M\eps^2,\ |\psi_1|+|\psi_2|\le M\eps,\\
  &\hspace{3cm} |\phi_1-\phi_2|\le M\eps^3,\ |\psi_1-\psi_2|\le M\eps^2\big\} \subset (H^1_{2\pi})^4.
\end{aligned}
\end{equation}
Observe that $(\phi,\psi)\in\cK_{M,\eps}$ implies $\big(u_1(t),\eps v_1(t),u_2(t),\eps v_2(t)\big) \in \cA_{R,\eps}$ with $R=M+B$ when $|\be|\le B$. Throughout this subsection $M,B>0$ are fixed, and all constants are uniform for $(\phi,\psi)\in\cK_{M,\eps}$ and $|\be|\le B$, for sufficiently small $\eps>0$. We need to control the terms on the right hand side of \eqref{eq:remainder-system}. First an elementary computation shows that we have in $\cK_{M,\eps}$:
\begin{equation}\label{eq:sing-u-est}
 \big|A_{\phi}\big| = O(\eps^4)\quad\text{and}\quad
 \big|A_\psi\big| = O(\eps^3)
\end{equation}

The metric factor $\ga_i-1$ appearing in $B_\phi$ and $B_\psi$ contributes an additional forcing of order $\eps^2$ by \eqref{eq:gamma-def}:
\begin{equation}\label{eq:metric-common}
 \gamma_i-1 = \frac12\kappa(\ubar)\eps^2+O(\eps^3).
\end{equation}
This contribution is of leading order in the equation for $\phi_1+\phi_2$ below and must be retained.

Lemma~\ref{lem:distance-remainder} and \eqref{eq:distance-corrections} give
\begin{equation}\label{eq:R-scaled-size}
 |\cR_{\psi i}|\le C\eps^3,\qquad  |\cR_{\phi i}|\le C\eps^5\qquad\text{for $i=1,2$.}
\end{equation}
Finally, \eqref{eq:F-gradient-structure}--\eqref{eq:F-gradient-sum} imply
\begin{equation}\label{eq:F-errors-size}
 |\eps\pa_{y_i}F|\le C\eps^2,
 \qquad
 |\pa_{x_i}F|\le C\eps,\qquad\text{for $i=1,2$.}
\end{equation}
and, more importantly,
\begin{equation}\label{eq:F-errors-difference}
 \eps|\pa_{y_1}F+\pa_{y_2}F|\le C\eps^3,
 \qquad
 |\pa_{x_1}F+\pa_{x_2}F|\le C\eps^2.
\end{equation}
The factors $\ga_i$ multiplying the distance corrections and the derivatives of $F$ in \eqref{eq:remainder-system} are retained. Their common correction satisfies $\ga_i-1=O(\eps^2)$, while $\ga_1-\ga_2=O(\eps^3)$; hence they preserve the orders below. 
Combining \eqref{eq:R-scaled-size}, \eqref{eq:gamma-difference}, and \eqref{eq:sing-u-est}--\eqref{eq:F-errors-difference}, the full right-hand side in
\[
 L^\eps(\phi,\psi)=\mathcal N^\eps(\phi,\psi,\be)
\]
has the natural sizes
\begin{equation}\label{eq:rhs-sizes}
 \cN^\eps_{\phi i}=O(\eps^2),
 \quad
 \cN^\eps_{\psi i}=O(\eps),
 \quad
 \cN^\eps_{\phi1}-\cN^\eps_{\phi2} = O(\eps^3),
 \quad
 \cN^\eps_{\psi1}-\cN^\eps_{\psi2} = O(\eps^2).
\end{equation}

%
\subsection{The projected system}
In order to solve the remainder system \eqref{eq:remainder-system} $L^\eps(\phi,\psi)=\cN^\eps(\phi,\psi,\be)$ we first project the right hand side onto the range of $L^\eps$ by setting
\begin{equation}\label{eq:c12}
 c_1 = \frac1{4\pi}\int_0^{2\pi}(\cN^\eps_{\phi2}-\cN^\eps_{\phi1})\,\dd t \quad\text{and}\quad
 c_2 = \frac1{4\pi}\int_0^{2\pi}(\cN^\eps_{\psi2}-\cN^\eps_{\psi1})\,\dd t,
\end{equation}
obtaining the projected system
\begin{equation}\label{eq:projected}
\left\{
\begin{aligned}
 L^\eps_{\phi1}(\phi,\psi) &= \cN^\eps_{\phi1}(\phi,\psi,\be) + c_1,\\
 L^\eps_{\phi2}(\phi,\psi) &= \cN^\eps_{\phi2}(\phi,\psi,\be) - c_1,\\
 L^\eps_{\psi1}(\phi,\psi) &= \cN^\eps_{\psi1}(\phi,\psi,\be) + c_2,\\
 L^\eps_{\psi2}(\phi,\psi) &= \cN^\eps_{\psi2}(\phi,\psi,\be) - c_2.
\end{aligned}
\right.
\end{equation}
We want to solve this system with $c_1=0=c_2$ for $\eps>0$ small. At first sight it seems that the problem is overdetermined since for the two scalar equations $c_i=0$ there is only one scalar unknown $\be$. However, the Hamiltonian character of the problem implies the following relation between $c_1$ and $c_2$.

\begin{Prop}\label{prop:c2-redundant}
Let $(\phi,\psi,\be)$ be a solution of \eqref{eq:projected}. Then
\begin{equation}\label{eq:energy-projected}
 c_1 I_x+\eps c_2 I_y=0,
\end{equation}
where
\begin{align*}
 I_x&=\int_0^{2\pi}
 \big(\pa_{x_1}\cH_\Phi(u_1,\eps v_1,u_2,\eps v_2)-\partial_{x_2}\cH_\Phi(u_1,\eps v_1,u_2,\eps v_2)\big)\,\dd t,\\
 I_y&=\int_0^{2\pi}
 \big(\pa_{y_1}\cH_\Phi(u_1,\eps v_1,u_2,\eps v_2)-\partial_{y_2}\cH_\Phi(u_1,\eps v_1,u_2,\eps v_2)\big)\,\dd t,
\end{align*}
with $\cH_\Phi$ from \eqref{eq:def-HPhi}. Moreover, for $(\phi,\psi)\in\cK_{M,\eps}$ and $|\be|\le B$, with $M,B>0$ fixed, the following estimate holds uniformly:
\begin{equation}\label{eq:Iy-nonzero}
 \eps I_y=-2+O(\eps^2)\qquad\text{as $\eps\to0$.}
\end{equation}
In particular $I_y\neq0$ for all sufficiently small $\eps$, and therefore
\[
 c_1=0\quad\Longrightarrow\quad c_2=0.
\]
\end{Prop}

\begin{proof}
Observe that $w_i(t):=\big(u_i(t),\eps v_i(t)\big)$ solves the following modified version of \eqref{HS-Fermi}:
\[
\left\{
\begin{aligned}
  \dot w_1 &= \eps J_{w_1} \nabla^g_{w_1}\cH_\Phi(w_1,w_2)+(c_1,\eps c_2)\\
  \dot w_2 &= -\eps J_{w_2} \nabla^g_{w_2}\cH_\Phi(w_1,w_2)-(c_1,\eps c_2)
\end{aligned}
\right.
\]
The Hamiltonian contributions to the derivative of the energy vanish because $g(\nabla^g\cH_\Phi,J\nabla^g\cH_\Phi)=0$ in each factor. Hence
\[
  \frac{d}{dt}\cH_\Phi(w_1,w_2) = c_1(\pa_{x_1}\cH_\Phi-\pa_{x_2}\cH_\Phi)(w_1,w_2) + \eps c_2(\pa_{y_1}\cH_\Phi-\pa_{y_2}\cH_\Phi)(w_1,w_2).
\]
The function $t\mapsto\cH_\Phi(w_1(t),w_2(t))$ is $2\pi$-periodic. because  $u_i(t+2\pi)=u_i(t)+1$ and $v_i(t+2\pi)=v_i(t)$ by \eqref{eq:periodicity} and $\Phi(x+1,y)=\Phi(x,y)$. Integrating the above identity over the period $[0,2\pi]$ of $\cH_\Phi(w_1,w_2)$ implies \eqref{eq:energy-projected}.

We now estimate $I_y$.  From \eqref{eq:def-HPhi} and Lemma~\ref{lem:distance-remainder} we deduce that
\[
 \pa_{y_1}\cH_\Phi(w_1,w_2)-\pa_{y_2}\cH_\Phi(w_1,w_2) = -\frac{2\eps(v_1-v_2)}  {\pi\bigl((u_1-u_2)^2+\eps^2(v_1-v_2)^2\bigr)}+O(\eps)
\]
uniformly on $\cK_{M,\eps}$.  Since 
\[
  v_1-v_2=2+\psi_1-\psi_2,\quad u_1-u_2=\phi_1-\phi_2
\]
and
\[
  \phi_1-\phi_2=O(\eps^3),\quad \psi_1-\psi_2=O(\eps^2),
\]
we have
\[
 \begin{aligned}
 \eps\big(\pa_{y_1}\cH_\Phi(w_1,w_2)-\pa_{y_2}\cH_\Phi(w_1,w_2)\big)
 &=-\frac{2\eps^2(2+\psi_1-\psi_2)}
 {\pi\bigl((\phi_1-\phi_2)^2+\eps^2(2+\psi_1-\psi_2)^2\bigr)}+O(\eps^2)\\
 &=-\frac{2}{\pi(2+\psi_1-\psi_2)}+O(\eps^2)
 =-\frac1\pi+O(\eps^2),
 \end{aligned}
\]
uniformly in $t$.  Integration over $[0,2\pi]$ yields \eqref{eq:Iy-nonzero}.  Notice in particular that $I_y$ is not zero; rather,
$I_y=-2/\eps+O(\eps)$.  If $c_1=0$, \eqref{eq:energy-projected} and \eqref{eq:Iy-nonzero} therefore imply $c_2=0$.
\end{proof}

\subsection{Conclusion of the proof}

In order to solve \eqref{eq:projected} with $c_1=0=c_2$ we pass to the equations for the sums and differences
\[
  S_\phi := \phi_1+\phi_2,\quad D_\phi := \phi_1-\phi_2,\quad S_\psi := \psi_1+\psi_2,\quad D_\psi := \psi_1-\psi_2,
\]
so that we can exploit the symmetries of $F$ as stated in Lemma~\ref{lem:F-identities} and its consequences \eqref{eq:F-gradient-structure}-\eqref{eq:Fxsum-centered}. We also consider the two scalar equations that $S_\phi$ and $S_\psi$ have mean value $0$, so that $(\phi,\psi)\in K^\perp$; see \eqref{eq:Kperp}. The equations involve three scalar variables $\be$, $c_1$, $c_2$. We set $c_1=0$ and keep $\be$ and $c_2$ as parameters. Proposition~\ref{prop:c2-redundant} then yields $c_2=0$ for a solution. The following rescaling of the sums, the differences, and $c_2$ will allow us to pass to the limit $\eps\to0$:
\[
  S_\phi =\eps^2\wt{S}_\phi,\quad D_\phi = \eps^3\wt{D}_\phi,\quad 
  S_\psi = \eps \wt{S}_\psi,\quad D_\psi = \eps^2\wt{D}_\psi,\quad c_2 = \eps^2 \wt{c}_2.
\]
Equivalently we have
\begin{equation}\label{eq:recover-phi-psi}
\left\{
\begin{aligned}
 \phi_1&=\frac{\eps^2\wt{S}_\phi+\eps^3\wt{D}_\phi}{2},&
 \phi_2&=\frac{\eps^2\wt{S}_\phi-\eps^3\wt{D}_\phi}{2},\\
 \psi_1&=\frac{\eps \wt{S}_\psi+\eps^2\wt{D}_\psi}{2},&
 \psi_2&=\frac{\eps \wt{S}_\psi-\eps^2\wt{D}_\psi}{2}.
\end{aligned}
\right.
\end{equation}

The substitutions in \eqref{eq:recover-phi-psi} now turn the system \eqref{eq:projected}, plus the two scalar equations corresponding to $(\phi,\psi)\in K^\perp$, into a system for the unknown
\[
  X := \big(\wt{S}_\phi,\wt{D}_\psi,\wt{c}_2, \wt{S}_\psi,\wt{D}_\phi,\be\big)
     \in H^1_{2\pi}\times H^1_{2\pi}\times\R\times H^1_{2\pi}\times H^1_{2\pi}\times\R =:  \mathbb X
\]
of the form
\begin{equation}\label{eq:scaled-exact-abstract}
 \mathscr A(X) = Y_0 + \mathscr E_\eps(X).
\end{equation}
Here the bounded linear operator
\[
 \mathscr A:\mathbb X\to\mathbb Y := L^2_{2\pi}\times L^2_{2\pi}\times\R\times L^2_{2\pi}\times L^2_{2\pi}\times\R
\]
is defined by
\begin{equation}\label{eq:A0}
\mathscr A(X) =
\left(
 \dot{\wt{S}}_\phi+\frac{\wt{D}_\psi}{2\pi},
 \dot{\wt{D}}_\psi-2\wt{c}_2,
 \int_0^{2\pi}\wt{S}_\phi,
 \dot{\wt{D}}_\phi-\frac{\ka(\ubar)}{2\pi}(2\be+\wt{S}_\psi),
 \dot{\wt{S}}_\psi+\frac{\wt{D}_\phi}{2\pi},
 \int_0^{2\pi}\wt{S}_\psi
\right).
\end{equation}
The element $Y_0\in \mathbb Y$ is defined as
\[
  Y_0 := \left(\frac{\ka(\ubar)}{2\pi} - 4c(\ubar,0),2\pa_x c(\ubar,0), 0, \frac{\ka_\nu(\ubar)}{6\pi}-2\pa_y c(\ubar,0), 4b(\ubar,0), 0\right) 
\]
where $b,c:\cS\to\R$ are the components of $A=D^2_{z_1,z_1} F:\cS\to\R^{2\times2}$ as in \eqref{def:beta_0}. In order to apply the contraction mapping principle to \eqref{eq:scaled-exact-abstract} we need to investigate when $\mathscr A$ is invertible, and to prove that $\mathscr A^{-1}\circ\mathscr E_\eps$ is a contraction.

\begin{Prop}\label{prop:A0}
The operator $\mathscr A$ is continuously invertible if and only if the operator 
\[
  \cJ_\Ga(h) = -h''-\ka h,
\]
has trivial kernel on $1$-periodic functions $h:\R\to\R$.
\end{Prop}

\begin{proof}
Observe that $\mathscr A=\begin{pmatrix}\mathscr A_1 &0\\0&\mathscr A_2\end{pmatrix}$ is block diagonal with respect to $(\wt S_\phi,\wt D_\psi,\wt c_2)$ and $(\wt S_\psi,\wt D_\phi,\be)$. For  $\mathscr A_1$ consider the system
\[
  \dot{\wt S}_\phi+\frac{\wt D_\psi}{2\pi} = f,\qquad \dot{\wt D}_\psi-2\wt c_2 = g,\qquad \int_0^{2\pi}\wt S_\phi = \al
\]
with $f,g\in L^2_{2\pi}$, $\al\in\R$. The second equation determines $\wt c_2\in\R$ uniquely because the periodicity of $\wt D_\psi$ implies that $2\wt c_2+g$ must have mean value $0$. This equation also determines ${\color{red}\wt D_\psi} \in H^1_{2\pi}$ up to a constant. This constant is determined in the same way by the first equation because $\wt S_\phi$ is periodic. Finally the first equation determines $\wt S_\phi \in H^1_{2\pi}$ up to a constant which is unique because the mean value of $\wt S_\phi$ is prescribed. Thus we have proved that 
\[
  \mathscr A_1: H^1_{2\pi}\times H^1_{2\pi}\times\R \to L^2_{2\pi}\times L^2_{2\pi}\times\R
\]
is bijective, hence continuously invertible by the closed graph theorem.

Concerning $\mathscr A_2$ we consider the system
\[
   \dot{\wt D}_\phi-\frac{\ka(\ubar)}{2\pi}(2\be+\wt S_\psi) = f,\quad 
   \dot{\wt S}_\psi+\frac{\wt D_\phi}{2\pi} = g,\quad 
   \int_0^{2\pi}\wt S_\psi = \al.
\]
Here $f,g\in L^2_{2\pi}$ and $\al\in\R$ are arbitrary. We pass to $1$-periodic functions by setting $x=t/2\pi$, $\wt f(x)=f(2\pi x)$, $\wt g(x)=g(2\pi x)$, and $Z(x) = 2\be+\wt S_\psi(2\pi x)$. Recall that $\ka\big(\ubar(t)\big) = \ka(t/2\pi) = \ka(x)$. The second equation now gives
\[
 \wt D_\phi(2\pi x) = 2\pi\wt g(x) - Z'(x).
\]
Substitution into the first equation yields, in the weak sense,
\begin{equation}\label{eq:Jacobi-inhom-augmented}
  \cJ_\Ga Z = -Z'' -\ka Z = 2\pi\bigl(\wt f-\wt g'\bigr) \qquad\hbox{in }H^{-1}_{\rm per}.
\end{equation}
Observe that $\cJ_\Ga:H^1_{\rm per}\to H^{-1}_{\rm per}$ is an isomorphism because $\ker \cJ_\Ga$ is trivial by assumption and $\cJ_\Ga$ is a self-adjoint elliptic operator. Therefore \eqref{eq:Jacobi-inhom-augmented} has a unique solution $Z\in H^1_{\rm per}$. The mean constraint determines $\be$ uniquely through
\begin{equation}\label{eq:beta}
 \be = \frac1{4\pi} \left(2\pi\int_0^1 Z(x)\,\dd x-\al\right).
\end{equation}
We reconstruct the unknowns by
\[
 \wt S_\psi(t)=z(t/2\pi)-2\be,\qquad
 \wt D_\phi(t)=2\pi g(t)-Z'(t/2\pi).
\]
The first function belongs to $H^1_{2\pi}$ and has integral $\al$ by
\eqref{eq:beta}. The second initially belongs to $L^2_{2\pi}$; however,
\eqref{eq:Jacobi-inhom-augmented} implies
\[
 \frac{\dd}{\dd x}\bigl(2\pi\wt g-Z'\bigr) = 2\pi\wt f+\ka Z\in L^2_{\rm per}.
\]
Consequently $\wt D_\phi\in H^1_{2\pi}$, and the reconstructed functions
satisfy both equations. Their uniqueness follows from the uniqueness of $Z$
and formula \eqref{eq:beta}. Thus $\mathscr A_2$ is bijective and has a bounded inverse.

Conversely, suppose that $0\ne h\in\ker\cJ_\Ga$. By elliptic regularity $h$ is smooth. Set
\[
 \be=\frac12\int_0^1h(x)\,\dd x,\qquad
 \wt S_\psi(t)=h(t/2\pi)-2\be,\qquad
 \wt D_\phi(t)=-h'(t/2\pi).
\]
This is a nonzero element of $\ker\mathscr A_2$, since $2\be+\wt S_\psi(2\pi x)=h(x)$, and hence $\mathscr A$ is not invertible.
This proves the equivalence.
\end{proof}

The nonlinearity $\mathscr E_\eps$ in \eqref{eq:scaled-exact-abstract} can be bounded as follows.

\begin{Lem}\label{lem:C1-remainder}
For every $M>0$ there holds:
\begin{equation}\label{eq:C1-small}
 \sup_{\|X\|_{\mathbb X}\le M}
 \left(
  \|\mathscr E_\eps(X)\|_{\mathbb Y} + \|D\mathscr E_\eps(X)\|_{\cL(\mathbb X,\mathbb Y)}\right)
     \longrightarrow0\qquad\text{as }\eps\to0.
\end{equation}
\end{Lem}

Before giving the proof of Lemma~\ref{lem:C1-remainder} we conclude the proof of the main theorem.\\

\begin{altproof}{Theorem~\ref{thm:main}}
By Proposition~\ref{prop:A0} we may rewrite equation~\eqref{eq:scaled-exact-abstract} as 
\[
 X=\mathscr A^{-1}(Y_0)+\mathscr A^{-1}\circ\mathscr E_\eps(X) =: \mathscr F_\eps(X).
\]
We seek a solution near $X_0:=\mathscr A^{-1}(Y_0)$. Set
$M:=1+\|X_0\|_{\mathbb X}$ and fix $0<r<1$.
Lemma~\ref{lem:C1-remainder} implies that, for sufficiently small $\eps>0$,
\[
  \sup_{X\in\overline B_r(X_0)}
  \|\mathscr A^{-1}\mathscr E_\eps(X)\|_{\mathbb X} \le r/2,
\]
and
\[
  \sup_{X\in\overline B_r(X_0)} \|\mathscr A^{-1}D\mathscr E_\eps(X)\|_{\cL(\mathbb X,\mathbb X)} \le 1/2.
\]
As a consequence $\mathscr F_\eps$ maps the closed ball $\overline B_r(X_0)$ into itself
and is a contraction there. This ball is complete, so the contraction mapping principle yields a unique fixed point 
\[
  X_\eps = \big(\wt S_\phi^\eps,\wt D_\psi^\eps,\wt c_2^\eps,\wt S_\psi^\eps,\wt D_\phi^\eps, \be_\eps\big)
\]
of $\mathscr F_\eps$, hence a solution of \eqref{eq:scaled-exact-abstract}. Clearly $\mathscr E_\eps(x)$, hence $\mathscr F_\eps(X)$ depend smoothly on $(\eps,X)$ in a neighborhood of these solutions. The bound $\|D_X\mathscr F_\eps\| \le 1/2$
implies that $I-D_X\mathscr F_\eps$ is invertible. The implicit function theorem
therefore gives smooth dependence of $X_\eps$ on $\eps\in(0,\eps_0)$. Moreover,
$\|X_\eps-X_0\|_{\mathbb X}
\le\|\mathscr A^{-1}\|\,\|\mathscr E_\eps(X_\eps)\|_{\mathbb Y}\to0$
by Lemma~\ref{lem:C1-remainder}. Now using the relations \eqref{eq:recover-phi-psi} we obtain a family of solutions $(\phi_1^\eps,\phi_2^\eps,c_2^\eps,\psi_1^\eps,\psi_2^\eps,\be_\eps)$ of the system~\eqref{eq:projected} with $c_1=0$. Proposition~\ref{prop:c2-redundant} implies $c_2^\eps=0$, hence $(\phi_1^\eps,\phi_2^\eps,\psi_1^\eps,\psi_2^\eps,\be_\eps)$ solves \eqref{eq:remainder-system}. Then
\[
\left\{
\begin{aligned}
 u_1^\eps &= \ubar+\phi_1^\eps, &&\qquad u_2^\eps = \ubar+\phi_2^\eps,\\
 v_1^\eps &= 1+\be_\eps\eps+\psi_1^\eps, &&\qquad  v_2^\eps = -1+\be_\eps\eps+\psi_2^\eps.
\end{aligned}
\right.
\]
solve \eqref{eq:scaled-system}, and
\[
  x_i(t) = u_i(t/\eps), \quad y_i(t) = \eps v_i(t/\eps),
\]
solve the original problem \eqref{HS} with $\ka_1=1=-\ka_2$. The equations \eqref{eq:period}--\eqref{eq:form y} are satisfied by construction, and the bounds \eqref{eq:bound1}--\eqref{eq:bound2} are satisfied by the relations~\eqref{eq:recover-phi-psi} and the fact that $X_\eps$ is bounded in $\mathbb X$, which implies that the $H^1_{2\pi}$ norms of $\wt S_\phi^\eps,\wt D_\psi^\eps,\wt S_\psi^\eps,\wt D_\phi^\eps$ are bounded. Finally, let $Z_0$ denote the function associated with the second block of $X_0$ as in the proof of Proposition~\ref{prop:linear}. Substituting the fourth and fifth components of $Y_0$ into
\eqref{eq:Jacobi-inhom-augmented} gives
\[
 -Z_0''-\ka Z_0
 =\frac{\ka_\nu}{3}-4\pi\pa_yc(\cdot,0)-8\pi\pa_xb(\cdot,0).
\]
This defines $Z_0$ as in \eqref{eq:Zsolution}. Since the sixth component of $Y_0$ is zero, \eqref{eq:beta} yields
$\be_0=\tfrac12\int_0^1Z_0(x)\,\dd x$.
The convergence $X_\eps\to X_0$ therefore proves
$\be_\eps\to\be_0$ with precisely the value in \eqref{eq:beta-limit-final}.
\end{altproof}

\begin{altproof}{Lemma~\ref{lem:C1-remainder}}
Throughout the proof $C$ denotes constants which may depend on $M$ but not on $\eps$ or $t$. We proceed in four steps: we first write $\mathscr E_\eps$ as a superposition (Nemytskii) operator, then show that the corresponding pointwise function is smooth up to $\eps=0$, then identify its value at $\eps=0$ by explicit expansions, and finally deduce \eqref{eq:C1-small}.

\medskip\noindent
\emph{Step 1: Pointwise form of the system.}
We want to write \eqref{eq:scaled-exact-abstract} as a pointwise differential equation $\cA_\eps(X)(t)=h(\eps,t,X(t))$. The definition of $h(\eps,t,\xi)\in\R$ proceeds in several steps, following the change of variables from $(u_1,u_2,v_1,v_2)\in (H^1_{2\pi})^4$ to $X\in\mathbb X$. For $t\in\R$ and $\xi=(s_\phi,d_\psi,\tau,s_\psi,d_\phi,\be)\in\R^6$ we set, in accordance with \eqref{eq:ansatz} and \eqref{eq:recover-phi-psi},
\begin{equation}\label{eq:pointwise-points}
\left\{
\begin{aligned}
 u_{1}&=\ubar(t)+\tfrac12\eps^2s_\phi+\tfrac12\eps^3d_\phi, &\qquad \eps v_1&=\eps+\eps^2\be+\tfrac12\eps^2s_\psi+\tfrac12\eps^3d_\psi,\\
 u_{2}&=\ubar(t)+\tfrac12\eps^2s_\phi-\tfrac12\eps^3d_\phi, & \eps v_2&=-\eps+\eps^2\be+\tfrac12\eps^2s_\psi-\tfrac12\eps^3d_\psi,
\end{aligned}
\right.
\end{equation}
and $z_i=(u_i,\eps v_i)$. Thus $\phi_1-\phi_2=\eps^3d_\phi$, $\psi_1-\psi_2=\eps^2d_\psi$, $\psi_1+\psi_2=\eps s_\psi$, and
\begin{equation}\label{eq:m-d}
\left\{
\begin{aligned}
  &m=\frac{z_1+z_2}2=\Big(\ubar(t)+\tfrac12\eps^2s_\phi,\ \eps^2\big(\be+\tfrac12s_\psi\big)\Big),\\
  &d=z_1-z_2=\eps w,\qquad w=\big(\eps^2d_\phi,\,2+\eps^2d_\psi\big).
\end{aligned}
\right.
\end{equation}
We denote by $N^\eps_{\phi i}(t,\xi)$, $N^\eps_{\psi i}(t,\xi)$ the right hand sides of \eqref{eq:remainder-system} evaluated at these points. Equivalently, by \eqref{eq:scaled-system},
\begin{equation}\label{eq:N-exact}
\left\{
\begin{aligned}
 N^\eps_{\phi1}(t,\xi) &= \ga_1\Big(\frac{\eps}{2\pi}\pa_{y_1}\log D_g-\eps\pa_{y_1}F\Big)(z_1,z_2)-\frac1{2\pi}+\frac{\eps^2d_\psi}{4\pi},\\
 N^\eps_{\phi2}(t,\xi) &= \ga_2\Big(-\frac{\eps}{2\pi}\pa_{y_2}\log D_g+\eps\pa_{y_2}F\Big)(z_1,z_2)-\frac1{2\pi}+\frac{\eps^2d_\psi}{4\pi},\\
 N^\eps_{\psi1}(t,\xi) &= \ga_1\Big(-\frac{1}{2\pi}\pa_{x_1}\log D_g+\pa_{x_1}F\Big)(z_1,z_2)+\frac{\eps d_\phi}{4\pi},\\
 N^\eps_{\psi2}(t,\xi) &= \ga_2\Big(\frac{1}{2\pi}\pa_{x_2}\log D_g-\pa_{x_2}F\Big)(z_1,z_2)+\frac{\eps d_\phi}{4\pi}.
\end{aligned}
\right.
\end{equation}
Define
\[
 g_1:=N^\eps_{\phi1}+N^\eps_{\phi2},\quad g_2:=N^\eps_{\psi1}-N^\eps_{\psi2},\quad
 g_4:=N^\eps_{\phi1}-N^\eps_{\phi2},\quad g_5:=N^\eps_{\psi1}+N^\eps_{\psi2},
\]
and
\begin{equation}\label{eq:def-h}
 h_1:=\eps^{-2}g_1,\quad h_2:=\eps^{-2}g_2,\quad h_4:=\eps^{-3}g_4-\frac{\ka(\ubar)}{2\pi}(2\be+s_\psi),\quad h_5:=\eps^{-1}g_5 .
\end{equation}
Taking the sum of the first two, the difference of the last two, the difference of the first two, and the sum of the last two equations in \eqref{eq:projected} with $c_1=0$, and using $c_2=\eps^2\wt c_2$, we obtain for $X\in\mathbb X$:
\[
\begin{aligned}
 &\dot{\wt S}_\phi(t)+\frac{\wt D_\psi(t)}{2\pi}=h_1\big(\eps,t,X(t)\big),&\quad
 &\dot{\wt D}_\psi(t)-2\wt c_2=h_2\big(\eps,t,X(t)\big),\\
 &\dot{\wt D}_\phi(t)-\frac{\ka(\ubar(t))}{2\pi}\big(2\be+\wt S_\psi(t)\big)=h_4\big(\eps,t,X(t)\big),&\quad
 &\dot{\wt S}_\psi(t)+\frac{\wt D_\phi(t)}{2\pi}=h_5\big(\eps,t,X(t)\big).
\end{aligned}
\]
Together with the two mean value conditions this is \eqref{eq:scaled-exact-abstract}, where the components of $\mathscr E_\eps$ are precisely
\begin{equation}\label{eq:def-E}
 \mathscr E_\eps(X)_j(t)=h_j\big(\eps,t,X(t)\big)-(Y_0)_j(t)\quad(j=1,2,4,5),\qquad \mathscr E_\eps(X)_3=\mathscr E_\eps(X)_6=0.
\end{equation}
Observe that $h_j$ does not depend on $\tau$, and that all functions are $2\pi$-periodic in $t$ because $\ubar(t+2\pi)=\ubar(t)+1$ and all geometric quantities are $1$-periodic in $x$.

\medskip\noindent
\emph{Step 2: Smoothness up to $\eps=0$.}
Fix $\rho>0$ and work on $\overline{B_\rho(0)}\subset\R^6$. We claim that there exists $\eps_0>0$ such that $g_1,g_2,g_4,g_5$ extend to $C^\infty$ functions on $(-\eps_0,\eps_0)\times\R\times\overline{B_\rho(0)}$ all of whose derivatives are uniformly bounded on compact $\eps$-intervals.
For $|\xi|\le\rho$ and $|\eps|$ small the points $z_1,z_2$ lie in $\cS$ and close to each other, so $\ga_i=1/\la(z_i)$ and $\nabla_{z_i}F(z_1,z_2)$ are smooth functions of $(\eps,t,\xi)$. The only singular ingredient is $\log D_g$. Since $D_g(z,z)=0$ and $\nabla D_g(z,z)=0$, Taylor's formula with integral remainder in the variable $d$ gives
\[
 D_g\big(m+\tfrac d2,m-\tfrac d2\big)=\big\langle G(m,d)d,d\big\rangle
\]
where
\[
  G(m,d)=\int_0^1(1-s)\,\nabla_d^2\big[D_g\big(m+\tfrac d2,m-\tfrac d2\big)\big]\big|_{sd}\,\dd s
\]
is smooth, even w.r.t.\ the variable $d$, and satisfies $G(m,0)=g(m)=\operatorname{diag}\big(\la^2(m),1\big)$. Inserting $d=\eps w$ from \eqref{eq:m-d} we get
\[
 D_g(z_1,z_2)=\eps^2q,\quad\text{with}
 \quad q(\eps,t,\xi):=\big\langle G(m,\eps w)w,w\big\rangle,\qquad q(0,t,\xi)=4,
\]
so $q$ is smooth and bounded below by $3$ for $|\eps|$ small. Similarly, writing $\nabla_{z_1}=\frac12\nabla_m+\nabla_d$, $\nabla_{z_2}=\frac12\nabla_m-\nabla_d$, and differentiating $\langle G(m,d)d,d\rangle$, every term contains at least one factor $d=\eps w$, so that
\[
  \nabla_{z_i}D_g(z_1,z_2)=\eps\,P_i(\eps,t,\xi),\qquad
  P_i(0,t,\xi)=(-1)^{i-1}2g(\ubar,0)\begin{pmatrix}0\\2\end{pmatrix}=(-1)^{i-1}\begin{pmatrix}0\\4\end{pmatrix},
\]
with $P_i$ smooth. Hence $\eps\pa_{y_i}\log D_g=(P_i)_y/q$ is smooth. Moreover $\pa_{x_i}\log D_g=\eps^{-1}(P_i)_x/q$, and since 
$(P_i)_x(0,t,\xi)=0$, Hadamard's lemma shows that $(P_i)_x=\eps\wt P_i$ with $\wt P_i$ smooth. Thus $\pa_{x_i}\log D_g$ is smooth 
as well, and the claim follows from \eqref{eq:N-exact}.

We shall use Hadamard's lemma in the following form: if $g$ is smooth on $(-\eps_0,\eps_0)\times\R\times\overline{B_\rho(0)}$ and $g(\eps,t,\xi)=\eps^k\ell(t,\xi)+O(\eps^{k+1})$ for some $k\ge1$, then $\pa_\eps^ig(0,\cdot)=0$ for $i<k$, so that
\[
 \eps^{-k}g(\eps,t,\xi)=\frac1{(k-1)!}\int_0^1(1-s)^{k-1}\pa_\eps^kg(s\eps,t,\xi)\,\dd s
\]
is smooth with all derivatives bounded, and equals $\ell(t,\xi)$ at $\eps=0$.

\medskip\noindent
\emph{Step 3: Expansions.}
We now show, uniformly for $t\in\R$ and $|\xi|\le\rho$, where $\ka,\ka_\nu$ are evaluated at $\ubar(t)$ and $b,c$ and their derivatives at $(\ubar(t),0)$:
\begin{equation}\label{eq:g-expansions}
\left\{
\begin{aligned}
 g_1&=\eps^2\Big(\frac{\ka}{2\pi}-4c\Big)+O(\eps^3), &\qquad g_2&=2\eps^2\pa_xc+O(\eps^3),\\
 g_4&=\eps^3\Big(\frac{\ka}{2\pi}(2\be+s_\psi)+\frac{\ka_\nu}{6\pi}-2\pa_yc\Big)+O(\eps^4), & g_5&=4\eps b+O(\eps^3).
\end{aligned}
\right.
\end{equation}
By Step~2 and Hadamard's lemma, \eqref{eq:g-expansions} implies that $h_1,h_2,h_4,h_5$ are smooth with bounded derivatives on $(-\eps_0,\eps_0)\times\R\times\overline{B_\rho(0)}$ and that
\begin{equation}\label{eq:h-at-zero}
 h_j(0,t,\xi)=(Y_0)_j(t)\qquad\text{for all $\xi$ and $j=1,2,4,5$;}
\end{equation}
note that in $h_4$ the term $\frac{\ka}{2\pi}(2\be+s_\psi)$ cancels exactly, so that the limit is independent of $\xi$.

For $|\xi|\le\rho$ the points \eqref{eq:pointwise-points} lie in $\cA_{R,\eps}$ for a suitable $R=R(\rho)$, so Lemma~\ref{lem:distance-remainder} and the estimates \eqref{eq:sing-u-est}--\eqref{eq:F-errors-difference} apply. We use the decomposition \eqref{eq:remainder-system}--\eqref{eq:distance-corrections} and collect the ingredients.

\emph{Metric factors.} Set $\rho_i:=\ga_i-1=O(\eps^2)$. By \eqref{eq:gamma}, \eqref{eq:gamma-difference} and $\psi_1+\psi_2=\eps s_\psi$:
\[
 \rho_1+\rho_2=\ka\eps^2+O(\eps^3),\qquad
 \rho_1-\rho_2=\eps^3\Big(\ka(2\be+s_\psi)+\frac{\ka_\nu}3\Big)+O(\eps^4).
\]

\emph{Euclidean terms.} The common denominator in \eqref{eq:AB} equals $\eps^2\big((2+\eps^2d_\psi)^2+\eps^4d_\phi^2\big)$. Therefore
\[
 e_\phi:=\frac{\eps^2}{\pi}\cdot\frac{2+\psi_1-\psi_2}{(\phi_1-\phi_2)^2+\eps^2(2+\psi_1-\psi_2)^2}
   =\frac1{\pi(2+\eps^2d_\psi)}+O(\eps^4)=\frac1{2\pi}-\frac{\eps^2d_\psi}{4\pi}+O(\eps^4),
\]
which gives the exact cancellation of the constant and the linear term in $A_\phi$:
\[
 A_\phi=e_\phi-\frac1{2\pi}+\frac{\eps^2d_\psi}{4\pi}=O(\eps^4),\qquad B_{\phi i}=\rho_ie_\phi=\frac{\rho_i}{2\pi}+O(\eps^4).
\]
Likewise
\[
 e_\psi:=\frac1\pi\cdot\frac{\phi_2-\phi_1}{(\phi_1-\phi_2)^2+\eps^2(2+\psi_1-\psi_2)^2}
   =-\frac{\eps d_\phi}{\pi\big((2+\eps^2d_\psi)^2+\eps^4d_\phi^2\big)}=-\frac{\eps d_\phi}{4\pi}+O(\eps^3),
\]
so that $A_\psi=e_\psi+\frac{\eps d_\phi}{4\pi}=O(\eps^3)$ and $B_{\psi i}=\rho_ie_\psi=O(\eps^3)$.

\emph{Distance corrections.} By \eqref{eq:R-scaled-size}: $\eps\ga_i\cR_{\phi i}=O(\eps^6)$ and $\ga_i\cR_{\psi i}=O(\eps^3)$.

\emph{Regular part.} By \eqref{eq:m-d} we have $d=\big(\eps^3d_\phi,\,2\eps+\eps^3d_\psi\big)$ and $m=(\ubar,0)+O(\eps^2)$, hence $A(m)=A(\ubar,0)+O(\eps^2)$ and similarly for $\pa_xA$, $\pa_yA$. Taking the difference, resp.\ the sum, of the two lines in \eqref{eq:F-gradient-structure}, and using \eqref{eq:Fysum-centered}--\eqref{eq:Fxsum-centered}, we obtain
\[
\begin{aligned}
 \pa_{y_1}F-\pa_{y_2}F&=2\big[A(m)d\big]_y+O(\eps^2)=4\eps c+O(\eps^2),\\
 \pa_{y_1}F+\pa_{y_2}F&=\tfrac12\langle\pa_yA(m)d,d\rangle+O(\eps^4)=2\eps^2\pa_yc+O(\eps^4),\\
 \pa_{x_1}F-\pa_{x_2}F&=2\big[A(m)d\big]_x+O(\eps^3)=4\eps b+O(\eps^3),\\
 \pa_{x_1}F+\pa_{x_2}F&=\tfrac12\langle\pa_xA(m)d,d\rangle+O(\eps^4)=2\eps^2\pa_xc+O(\eps^4),
\end{aligned}
\]
and $|\nabla_{z_i}F|=O(\eps)$.

\emph{Metric weights.} For any $f_1,f_2$ there holds
\begin{equation}\label{eq:gamma-weights}
 \ga_1f_1\pm\ga_2f_2=\Big(1+\frac{\rho_1+\rho_2}2\Big)(f_1\pm f_2)+\frac{\rho_1-\rho_2}2\,(f_1\mp f_2).
\end{equation}
Since $\rho_1+\rho_2=O(\eps^2)$ and $\rho_1-\rho_2=O(\eps^3)$, this identity keeps track of both the common metric correction and the smaller difference of the metric factors. We apply it below to the sum and difference of each pair of derivatives of $F$; for the distance corrections, the individual bounds already suffice.

We can now verify \eqref{eq:g-expansions}.

\emph{Component $1$.} Adding the first two equations of \eqref{eq:remainder-system} (note the opposite signs of the $F$-terms, which reflect the opposite circulations in \eqref{eq:scaled-system}),
\[
 g_1=2A_\phi+(\rho_1+\rho_2)e_\phi+\eps\big(\ga_1\cR_{\phi1}+\ga_2\cR_{\phi2}\big)-\eps\big(\ga_1\pa_{y_1}F-\ga_2\pa_{y_2}F\big).
\]
By \eqref{eq:gamma-weights}, $\ga_1\pa_{y_1}F-\ga_2\pa_{y_2}F=4\eps c+O(\eps^2)$. Hence
$g_1=O(\eps^4)+\frac{\ka\eps^2}{2\pi}+O(\eps^3)+O(\eps^6)-4c\eps^2+O(\eps^3)$, which is the first formula in \eqref{eq:g-expansions}.

\emph{Component $2$.} Subtracting the last two equations, $A_\psi$ cancels and
\[
 g_2=(\rho_1-\rho_2)e_\psi+\big(\ga_1\cR_{\psi1}-\ga_2\cR_{\psi2}\big)+\big(\ga_1\pa_{x_1}F+\ga_2\pa_{x_2}F\big).
\]
The first term is $O(\eps^4)$, the second is $O(\eps^3)$, and by \eqref{eq:gamma-weights} the third equals $2\eps^2\pa_xc+O(\eps^4)$.

\emph{Component $4$.} Subtracting the first two equations, $A_\phi$ cancels and
\[
 g_4=(\rho_1-\rho_2)e_\phi+\eps\big(\ga_1\cR_{\phi1}-\ga_2\cR_{\phi2}\big)-\eps\big(\ga_1\pa_{y_1}F+\ga_2\pa_{y_2}F\big).
\]
Here $(\rho_1-\rho_2)e_\phi=\frac{\rho_1-\rho_2}{2\pi}+O(\eps^5)=\eps^3\big(\frac{\ka}{2\pi}(2\be+s_\psi)+\frac{\ka_\nu}{6\pi}\big)+O(\eps^4)$, the second term is $O(\eps^6)$, and by \eqref{eq:gamma-weights} $\eps(\ga_1\pa_{y_1}F+\ga_2\pa_{y_2}F)=2\eps^3\pa_yc+O(\eps^5)$. This is the third formula in \eqref{eq:g-expansions}. Observe that the leading terms of order $\eps^2$, which are present in $\cN^\eps_{\phi1}$ and $\cN^\eps_{\phi2}$ separately, cancel in the difference; this is the reason why $D_\phi$ can be scaled by $\eps^3$.

\emph{Component $5$.} Adding the last two equations,
\[
 g_5=2A_\psi+(\rho_1+\rho_2)e_\psi+\big(\ga_1\cR_{\psi1}+\ga_2\cR_{\psi2}\big)+\big(\ga_1\pa_{x_1}F-\ga_2\pa_{x_2}F\big).
\]
The first three terms are $O(\eps^3)$, and by \eqref{eq:gamma-weights} the last one equals $4\eps b+O(\eps^3)$.

\emph{Components $3$ and $6$.} These are the mean value conditions $\int_0^{2\pi}\wt S_\phi=0=\int_0^{2\pi}\wt S_\psi$ expressing $(\phi,\psi)\in K^\perp$; they are linear and exact, so $\mathscr E_\eps(X)_3=\mathscr E_\eps(X)_6=0$.

\medskip\noindent
\emph{Step 4: Conclusion.}
By \eqref{eq:def-E} and \eqref{eq:h-at-zero} we have $\mathscr E_\eps(X)_j(t)=e_j\big(\eps,t,X(t)\big)$ with
\[
 e_j(\eps,t,\xi):=h_j(\eps,t,\xi)-h_j(0,t,\xi)=\eps\int_0^1\pa_\eps h_j(s\eps,t,\xi)\,\dd s.
\]
Since $h_j$ is smooth with bounded derivatives, there exists $C=C(\rho)$ such that
\begin{equation}\label{eq:pointwise-e}
 |e_j(\eps,t,\xi)|+|\pa_\xi e_j(\eps,t,\xi)|\le C\eps\qquad\text{for $0<\eps<\eps_0$, $t\in\R$, $|\xi|\le\rho$.}
\end{equation}
Now let $C_S$ be a constant for the continuous embedding $\mathbb X\hookrightarrow C^0_{2\pi}(\R^6)$, identifying the scalar components with constant functions, and choose $\rho:=1+C_SM$, so that $\|X\|_{\mathbb X}\le M$ implies $|X(t)|\le\rho$ for all $t$. The superposition operator $X\mapsto e_j(\eps,\cdot,X(\cdot))$ is of class $C^1$ from $\mathbb X$ into $L^2_{2\pi}$, with derivative
\[
 D\mathscr E_\eps(X)_j[\de X](t)=\pa_\xi e_j\big(\eps,t,X(t)\big)\,\de X(t),
\]
because $\mathbb X$ embeds continuously into the space of bounded continuous functions with values in $\R^6$ (the components $\wt c_2,\be$ being constants), and $\pa_\xi e_j$ is uniformly continuous on compact sets. Therefore \eqref{eq:pointwise-e} yields
\[
 \|\mathscr E_\eps(X)\|_{\mathbb Y}\le C\eps,\qquad
 \|D\mathscr E_\eps(X)[\de X]\|_{\mathbb Y}\le C\eps\sup_t|\de X(t)|\le CC_S\,\eps\,\|\de X\|_{\mathbb X},
\]
for all $\|X\|_{\mathbb X}\le M$ and $0<\eps<\eps_0$, which proves \eqref{eq:C1-small}.

Finally we note that, since the $h_j$ are smooth in $(\eps,\xi)$ up to $\eps=0$, the map $(\eps,X)\mapsto\mathscr E_\eps(X)$ is smooth from $(-\eps_0,\eps_0)\times B_M(0)$ into $\mathbb Y$, with $\mathscr E_0\equiv0$. This justifies the smooth dependence of $X_\eps$ on $\eps$ used in the proof of Theorem~\ref{thm:main}. Alternatively, the implicit function theorem applied to $(\eps,X)\mapsto\mathscr A(X)-Y_0-\mathscr E_\eps(X)$ at $(0,X_0)$, whose $X$-derivative is the isomorphism $\mathscr A$, yields the same branch $X_\eps$.
\end{altproof}

\section*{Acknowledgements}
The authors acknowledge the use of ChatGPT (OpenAI) to assist with selected calculations and with revising the manuscript. The main mathematical ideas, constructions and proof strategies are due to the authors. The authors have independently verified all mathematical arguments and calculations and take full responsibility for the content of the paper.


\begin{thebibliography}{99}

\bibitem{aref-etal:2002} H. Aref, P. K. Newton, M. A. Strember, T. Tokieda, D. Vainchtein:
{\em Vortex crystals.} Adv. Appl. Mech. {\bf 39} (2002), 1--79.

\bibitem{aref:2008} H. Aref:
{\em Point vortex dynamics: A classical mathematics playground.}
{\color{red}J. Math. Phys. {\bf 48} (2007), 065401.}

\bibitem{Avramidi:2015} I. Avramidi:
{\em Heat Kernel Methods and its Applications.}
Springer, Cham 2015.

\bibitem{bartsch-dai:2016} T. Bartsch, Q. Dai:
{\em Periodic solutions of the $N$-vortex Hamiltonian system in planar domains.}
J. Diff. Equ. {\bf 260} (2016), 2275--2295

\bibitem{bartsch-gebhard:2017} T. Bartsch, B. Gebhard:
{\em Global continua of periodic solutions of singular first-order Hamiltonian systems of $N$-vortex type.}
Math. Ann. {\bf 369} (2017), 627-651.

\bibitem{bartsch-sacchet:2018} T. Bartsch, M. Sacchet:
{\em Periodic solutions with prescribed minimal period of vortex type problems in domains.}
Nonlinearity {\bf 31} (2018), 2156--2172.

\bibitem{Koiller-Boatto:2009} S. Boatto, J. Koiller:
{\em Vortex pairs on surfaces.}
AIP Conf. Proc. {\bf 1130} (2009), 77--88.

\bibitem{Boatto-Koiller:2015} S. Boatto, J. Koiller: 
{\em Vortices on closed surfaces.} 
Geometry, mechanics, and dynamics, 185–237. Fields Inst. Commun., 73, Springer, New York, 2015.

\bibitem{Cao-Li-Wang:2025} D. Cao, S. Li, G. Wang:
{\em Desingularization of vortices for the incompressible Euler equation on a sphere.}
Calc. Var. Partial Differential quations {\bf 64} (2025), Paper No. 279, 26 pp.

\bibitem{Castro-Cordoba-Gomez:2019} A. Castro, D. C\'ordoba, J. G\'omez-Serrano:
{\em Uniformly rotating smooth solutions for the incompressible 2D Euler equations.}
Arch. Ration. Mech. Anal. {\bf 231} (2019), 719-785.

\bibitem{DdPMW:2020} J. Davila, M. del Pino, M.Muss, J. Wei:
{\em Gluing methods for vortex dynamics in Euler flows.}
Arch. Ration. Mech. Anal. {\bf 235} (2020), 1467-1530.

\bibitem{Dritschel-Boatto:2015} D. G. Dritschel, S. Boatto:
{\em The motion of point vortices on closed surfaces.}
Proc. R. Soc. A {\bf 471} (2015), 20140890.

\bibitem{Drivas-Glukhovskiy-Khesin:2024} T. D. Drivas, D. Glukhovskiy, B. Khesin:
{\em Singular vortex pairs follow magnetic geodesics.}
Int. Math. Res. Not. IMRN {\bf 2024} (2024), no.~14, 10880--10894.

\bibitem{flucher-gustafsson:1997} M. Flucher, B. Gustafsson:
{\em Vortex motion in two-dimensional hydromechanics.} 
Preprint in TRITA-MAT-1997-MA-02.

\bibitem{Gelantalis-Sternberg:2012} M. Gelantalis, P. Sternberg:
{\em Rotating $2N$-vortex solutions to the Gross-Pitaevskii equation on $S^2$.}
J. Math. Phys. {\bf 53} (2012), 083701, 24 pp.

\bibitem{Gromoll-Klingenberg-Meyer:1975} D. Gromoll, W. Klingenberg, W. Meyer: 
{\em Riemannsche Geometrie im Gro\ss en.}
Lecture Notes in Mathematics 55, Springer-Verlag, Berlin, 1975.

\bibitem{GrottaRagazzo-Gustafsson-Koiller:2024}
C. Grotta-Ragazzo, B. Gustafsson, J. Koiller:
{\em On the interplay between vortices and harmonic flows: Hodge decomposition of Euler's equations in 2d.}
Regul. Chaotic Dyn. {\bf 29} (2024), 241--303.

\bibitem{Gustafsson:2022} B. Gustafsson:
{\em Vortex pairs and dipoles on closed surfaces.}
J. Nonlinear Sci. {\bf 32} (2022), Paper No.~62.

\bibitem{Helmholtz:1858} H. Helmholtz:
{\em \"{U}ber Integrale der hydrodynamischen Gleichungen, welche den
Wirbelbewegungen entsprechen.}
J. reine angew. Math. {\bf 55} (1858), 25--55.

\bibitem{Hmidi-Xue-Xue:2026} T. Hmidi, L. Xue and Z. Xue:
{\em Doubly connected V-states in geophysical models: a general framework.}
Calc. Var Partial Differential Equaations {\bf 65} (2026), Paper No. 137, 68 pp.

\bibitem{Kimura:1999} Y. Kimura:
{\em Vortex motion on surfaces with constant curvature.}
Proc. R. Soc. Lond. A {\bf 455} (1999), 245--259.

\bibitem{Kirchhoff:1876} G. Kirchhoff:
{\em Vorlesungen \"{u}ber mathematische Physik.}
Teubner, Leipzig 1876.

\bibitem{lin:1941a} C. C. Lin:
{\em On the motion of vortices in 2D I. Existence of the Kirchhoff-Routh function.}
Proc. Nat. Acad. Sc. {\bf 27} (1941), 570--575.

\bibitem{lin:1941b} C. C. Lin:
{\em On the motion of vortices in 2D II. Some further properties on the Kirchhoff-Routh function.}
Proc. Nat. Acad. Sc. {\bf 27} (1941), 575--577.

\bibitem{Long:2002} Y. Long:
 Index Theory for Symplectic Paths with Applications. Birkh\"auser Verlag, Basel 2002.

\bibitem{majda-bertozzi:2001} A. J. Majda, A. L. Bertozzi:
{\em Vorticity and Incompressible Flow.}
Cambridge University Press 2001.

\bibitem{marchioro-pulvirenti:1994} C. Marchioro, M. Pulvirenti:
{\em Mathematical Theory of Incompressible Nonviscous Fluids.}
Applied Mathematical Sciences {\bf 96}, Springer-Verlag, New York 1994.

\bibitem{newton:2001} P. K. Newton:
{\em The $N$-vortex problem.}
Springer-Verlag, Berlin 2001.

\bibitem{Rodrigues-Castilho-Koiller:2018} A. Regis Rodrigues, C. Castilho, J. Koiller:
{\em Vortex pairs on a triaxial ellipsoid and Kimura's conjecture.}
J. Geom. Mech. {\bf 10} (2018), 189--208.

\bibitem{routh:1881} E. J. Routh: {\em Some applications of conjugate functions.}
Proc. London Mat. Soc. {\bf 12} (1881), 73--89.

\bibitem{saffman:1992} P. G. Saffman: {\em Vortex dynamics.}
Cambridge University Press, Cambridge 1992.

\bibitem{Sakajo-Shimizu:2016} T. Sakajo, Y. Shimizu:
{\em Point vortex interactions on a toroidal surface.}
Proc. R. Soc. A {\bf 472} (2016), 20160271.

\bibitem{wang:2021} Q. Wang: 
{\em The $N$-vortex problem on a Riemann sphere.}
Commun. Math. Phys. {\bf 385} (2021), 565-593.

\end{thebibliography}
\end{document}